\documentclass[11pt,leqno]{article}

\usepackage{amsmath,amsthm}
\usepackage {latexsym}
\usepackage{amssymb}

\newcommand{\les}{\lesssim}
\newcommand{\bea}{\begin{eqnarray}}\newcommand{\eea}{\end{eqnarray}}
\newcommand{\beq}{\begin{equation}}\newcommand{\ee}{\end{equation}}

\renewcommand{\b}{\beta}

\newtheorem{theorem}{Theorem}[section]\newtheorem{lemma}[theorem]{Lemma}
\newtheorem{cor}[theorem]{Corollary}\newtheorem{prop}[theorem]{Proposition}

\theoremstyle{remark}

\def\ve{\varepsilon}\def\a{\alpha}\def\Si{\Sigma}
\def\bm{\left( \begin{array}{cc}}
\def\endm{\end{array}\right)}\newcommand{\eq}{\end{equation}}
\def\a{\alpha}\def\b{\beta}
\def\Box{\square}\def\pa{\partial}\def\pab{\bar\pa}

\def \rectangle#1#2{\hbox{\vrule\vbox to #2 {\hrule\hbox to #1{\hfil}\vfil\hrule}\vrule}}
\def\Lb{\underline{L}}
\def\a{\alpha}\def\b{\beta}
\def\Box{\square}\def\Boxr{\widetilde{\square}}\def\pa{\partial}
\def\pab{\bar\pa}\def\Lb{\underline{L}}

\def\pam{(\pa_{t\!}-_{\!}\pa_r)}\def\pap{(\pa_{t\!}+_{\!}\pa_r)}

\def\u#1{\underline{#1\!}\,}
\def\nup{{\nu^{\,\prime}}}

\numberwithin{equation}{section}

\begin{document}

\title {Global solutions of quasilinear wave equations}
\author {Hans Lindblad \\
University of California at San Diego} \maketitle

\section{Introduction}
We show that the Cauchy problem in $\bold{R}^{1+3}$:
\beq\label{eq:cauchy}
\Boxr_{g(\phi)}\phi =0,\qquad\qquad
\phi\big|_{t=0}=\phi_0,\quad \pa_t\phi\big|_{t=0}=\phi_1
\eq
has a global solution for all $t\geq 0$ if initial data
are sufficiently small.
Here the curved wave operator is
$
\Boxr_g=g^{\alpha\beta}\pa_\alpha\pa_\beta,
$
where we used the convention that repeated upper and lower indices are summed over $\alpha,\beta=0,1,2,3$, and $\pa_0=\pa/\pa t$, $\pa_i=\pa/\pa x^i$, $i=1,2,3$.
We assume that $g^{\alpha\beta}(\phi)$ are smooth functions of $\phi$
such that $g^{\alpha\beta}(0)\!=m^{\alpha\beta}$, where $m^{00}\!=\!-1$,
$m^{11}\!=m^{22}\!=m^{33}\!=\!1$ and $m^{\alpha\beta}\!=0$, if $\alpha\!\neq\!\beta$.
The result holds for vector valued $\phi$, in particular for the principal part of Einstein's equations; $\phi^{\alpha\beta}\!=g^{\alpha\beta}\!-m^{\alpha\beta}\!$.

This result was conjectured in \cite{L2} where it was also shown
in the spherically symmetric case for
\beq\label{eq:model}
-\pa_t^2 \phi+c(\phi)^2 \triangle_x\phi=0,\qquad\qquad\text{where}
\qquad c(0)=1.
\eq
In \cite{L2} there was also a heuristic argument for why the conjecture should be true in general:
Consider
\begin{equation}\label{eq:Quas}
\Box\,\phi = a_{\bold{\alpha\beta}}\, \pa^{\bf \alpha}\phi\,\pa^\beta\phi + {\text{cubic terms}},
\end{equation}
where here we used multiindex notation and the sum is over $0\!\leq\! |\alpha|\!\leq\! |\beta|\!\leq\! 2$ and $a_{\alpha\beta}$
are constants.
If we neglect derivatives tangential to the outgoing Minkowski light cones
and cubic terms, that are known to decay faster,
we get the {asymptotic equation} for $\Phi=r\phi$, introduced by H\"ormander
\cite{H1, H2, H3}:
\beq\label{eq:asymp}
\pap\pam\Phi \sim\\{r}^{-1} \! A_{mn}
\pam^m\Phi \,\, \pam^n \Phi,\qquad
A_{mn}=\!
\frac{1}{4}\! \sum_{|\alpha|=m|,\, |\beta|=n\!\!\!\!\!\!\!\!\!\!\!\!\!\!\!\!\!\!\!\!\!\!\!\!}
 a_{\alpha\beta}
\hat\omega^\alpha \hat \omega^\beta,\,\,\,\,\,\hat
\omega\!=\!(\!-\!1,\omega).
\eq
Here we have introduced polar coordinates $x=r\omega$, $\omega\in\bold{S}^2$.
The classical {\it null condition} introduced by Klainerman \cite{K1} is that $A_{nm} \equiv 0$ under which Klainerman \cite{K2} and Christodoulou \cite{C} proved global existence. In \cite{L2} it was observed that the asymptotic equation corresponding to \eqref{eq:cauchy} has global solution\footnote{In \cite{L-R1} we in general say that \eqref{eq:Quas} satisfy the {\it weak null condition}, if \eqref{eq:asymp} has global solution with some decay.},
contrary to other cases like $\Box \phi =\phi_t \triangle \phi$
or $\Box \phi= \phi_t^2$, where solutions are known to blow up for all small data, see John \cite{J1,J2}.
However, unlike for the classical null condition, the solution of \eqref{eq:cauchy} do not behave asymptotically like a solution of a free linear wave equation.

The method of proof of \cite{L2} is integration along characteristics
so it does not directly generalize to the non-symmetric case.
However, as observed in \cite{L1}, the method of integration along characteristics can still be used to obtain
sharp decay estimates assuming weaker decay estimates that can be obtained from
energy estimates for vector fields applied to the solution.
This then has to be combined with some refined energy estimates that take into account that the characteristic surfaces curve asymptotically,
since the solution do not decay as much as a solution of a free linear wave equation.

 Recently Alinhac \cite{A2}
 generalized the result in \cite{L2} to general data for the special
 case \eqref{eq:model}\footnote{As mentioned in \cite{A3} the method
 of \cite{A2} use the special structure of \eqref{eq:model}
 and is unlikely to work in general.}.
 \cite{A2} combines ideas from \cite{L1,L2} of how to obtain decay
 estimates with ideas from \cite{A1} for energy estimates with
 weights. Because the asymptotic behavior is different from that of solutions
 to a free linear wave equation, \cite{A2} constructs vector fields
 adapted to the characteristic surfaces at infinity,
 which in spirit is similar to the work of Christodoulou-Klainerman \cite{C-K}.
 Since these depend on the solution itself, commuting the vector fields with
 the wave operator leads to a loss of regularity
 so it has to be combined with a smoothing procedure, which leads to long
 schematic commutator estimates.

 There is however no need to construct vector fields adapted to the
 geometry at infinity. In fact we just use the vector
 fields for the Minkowski space time. In \cite{L-R3},
 for Einstein's equations, we
 also got away with just using the regular vector fields, but only because
 we got additional control from the wave coordinate condition.
 The observations here will hopefully will lead to a proof
 that uses less of the special structure and applies to a more general
 class of equations, which is useful in applications.

 As mentioned above the proof involves obtaining sharp decay estimates for low
 derivatives
 just assuming a weak decay estimate that later will be obtained from energy
 estimates for higher derivatives.
 The sharp decay estimates uses integration along characteristics as
 in \cite{L1,L2, L-R2, L-R3}. We adopt the energy method
 with weights of \cite{A2},
 depending on the solution of an approximate eikonal equation.
  This is a much
 easier substitute for energies on characteristic surfaces as I originally
 planned to use.
 The construction of vector fields adapted to the asymptotic behavior
 of the characteristic surfaces of \cite{A2} is
 avoided by considering a family of energy and decay estimates
 for the vector fields of flat Minkowski space time, with
 different decays for different types of derivatives.
 We prove the following:
 \begin{theorem}\label{mytheorem}  Suppose that $\phi_0$ and $\phi$ are smooth
 functions such that $\phi_0(x)=\phi_1(x)=0$,
 when $|x|\geq 1$, and let $N\geq 14$. Then there is a constant
 $\varepsilon_0>0$, such that if
 $\sum_{|\alpha|\leq N}
\big( \|\pa \pa^\alpha \phi_0\|_{L^2}+\|\pa^\alpha \phi_1\|_{L^2}\big)=\varepsilon\leq \varepsilon_0$,
 the Cauchy problem \eqref{eq:cauchy} has a global solution $\phi$
 for all $t\geq 0$.
The solution satisfies the decay estimates in Proposition \ref{strongdecay} with $\nu\!=\!1\!-\!c\varepsilon$, for some $c\!>\!0$, and the energy estimates
in Proposition \ref{energybound}.
 \end{theorem}

We remark that the result is true also for systems
$\phi=(\phi_1,\cdots,\phi_M)$, in particular the principal
quasilinear part of Einstein's equations. We also remark that the
assumptions on compact support is not needed and we can include
decaying data using energy norms with weights as in
\cite{L-R2,L-R3}.

Let us now give the strategy of the proof and the main ideas.
The proof involves getting sharp decay estimates for low derivatives
assuming weak decay estimates, and energy estimates for high derivatives assuming sharp decay estimates for low derivatives. The weak decay estimates
can then be obtained from energy estimates using
a bootstrap or continuity argument that we describe below.
Let
\beq
E_N(t)=\sum_{|I|\leq N}\int |\pa Z^I\phi(t,x)|^2\, dx,
\eq
where $Z^I$ is a product of $|I|$ of the vector fields,
$\Omega_{\alpha\beta}=x_\alpha \pa_\beta-x_\beta\pa_\alpha$, $S=x^\alpha\pa_\alpha$  that span the tangent space of the forward light cone and have good commutators with the wave operator, and $\pa_\alpha$.
(Here $x_i=x^i$, $i\geq 1$, $x_0=-x^0=-t$.)
In view of local existence results it suffices to give a bound for
$E_N(t)$, which will be obtained through a {\bf continuity argument}, see section \ref{proof}.
 Fix $0\!<\!\delta\!<\!1$. Assuming that
\beq\label{eq:energyassumption}
E_N(t)\leq 16 N \varepsilon^2 (1+t)^\delta,
\eq
for $0\leq t\leq T$, which holds for $T=0$, we will show that this bound,
implies the same bound with $16$ replaced by $8$ if $\varepsilon$ is
sufficient small (independently of $T$).
Using Klainerman-Sobolev inequality and the assumption of compactly supported data, see sections \ref{klainermansobolev}, \ref{proof},  this gives
{\bf weak decay estimates}:
\beq\label{eq:weakdecayintro}
|Z^I\phi(t,x)|
\leq c_0\varepsilon (1+t)^{-\nu}, \qquad
\nu>0, \quad
|I|\leq N-2.
\eq
These weak decay estimates imply
the sharp decay estimates in Proposition \ref{strongdecay},
as well as the estimates for the approximate radial characteristic surfaces
in Proposition \ref{eikonalone}, Lemma \ref{eikonaltwo} and Lemma \ref{eikonalthree}.
These sharp decay estimates for low derivatives are sufficient for the energy estimate in Proposition \ref{energybound} to hold
and we therefore get back a stronger energy estimate if $\varepsilon>0$ is
sufficient small:
\beq
E_N(t)\leq 8N\varepsilon^2(1+t)^{C_{0,N}\varepsilon}.
\eq

We now give the main ideas for the {\bf sharp decay estimates}.
We will try to mimic the integration along characteristic
that was done in the radial case in \cite{L2}, by expressing the wave operator in spherical coordinates and a null-frame, using the weak decay estimates to control the angular derivatives.
The discussion below will be a bit technical, but it is useful to get a feeling for how the different kind of terms are dealt with since the
structure of the argument is the same also for the energy estimates.

In  section \ref{nullframe} we express the inverse of the metric in terms of a nullframe:
\beq\label{eq:metricnullframeintro}
g^{\alpha\beta}=-\frac{1}{2}\big(\, L_1^\alpha \underline{L}^\beta+
 \underline{L}^\alpha L_1^{\beta}\big)+\gamma^{\alpha\beta},
\eq
where
\begin{align}
L_1^\alpha&=L^\alpha-H^{\u L\u L}\u L^\alpha -2H^{L\u L} L^\alpha
-2H^{\u L A} A^\alpha,\\
 \gamma^{\alpha\beta}&=\delta^{AB}A^\alpha B^\beta +H^{LL}L^\alpha L^\beta+H^{AL} A^\alpha L^\beta+ H^{AL} L^\alpha A^\beta +H^{AB} A^\alpha B^\beta.
\end{align}
Here $H^{\u L\u L}$ etc. are the components of $H^{\alpha\beta}=g^{\alpha\beta}-m^{\alpha\beta}$ in the Minkowski null frame
\beq
\underline{L}=(1,-\omega),\qquad L=(1,\omega), \qquad A,B\in \bold{S}^2,\qquad
\delta_{\alpha\beta}A^\alpha B^\beta=\delta_{AB}.
\eq

In section \ref{vectorfield} we use
\eqref{eq:metricnullframeintro} to decompose the wave operator:
 \beq \label{eq:waveframeintro}
\Big| \Big(2 L_1^{\alpha}\pa_\alpha
+\frac{\ell}{r}\Big)(r\pa_q\phi)
-r\Boxr_g \phi\Big|\les\sum_{1\leq |k|\leq
2\!\!\!\!\!\!\!\!\!\!} r^{|k|-1}\,|\pab^{\,k}\phi|,
\qquad\qquad\pa_q=\tfrac{1}{2}(\pa_r\!-\!\pa_t),
 \eq
where $\ell=\delta_{AB} H^{AB}\!
+ 4 H^{\u L L}\!-2 H^{\u L\u L}$, and $\overline{\pa}\in \{L,A,B\}$ are derivatives tangential to
the outgoing Minkowski light cones, that can be estimated in terms of the vector fields:
\beq\label{eq:Zpaq}
(1+|\,t-r|)|\pa\phi|+(1+t+r)|\overline{\pa}\phi|\les \sum_{|I|=1} |Z^I\phi|.
\eq
Note that when $|t-r|>t/2$ this together with \eqref{eq:weakdecayintro} gives
the sufficient $\varepsilon t^{-1}$ decay for all derivatives but when
$|t-r|$ is close to the light cone we are missing one derivative perpendicular to the light cone.

In section \ref{strongdecaywave} we integrate \eqref{eq:waveframeintro} along the flow lines of the vector field
$2L_1^\alpha\pa_\alpha$, from $|t-r|=t/2$,
to also get an estimate for a derivative perpendicular
to the outgoing light cones $r\pa_q \phi$ which yields
\begin{multline}\label{eq:strongdecaywaveeqintro}
(1+t+r) |\pa \phi(t,x)| \leq  C\sup_{0\leq \tau\leq \,t}
\sum_{|I|\leq 1}
\|Z^I\! \phi(\tau,\cdot)\|_{L^\infty}\\
+ C\int_{0}^t \Big((1+\tau)\| \Boxr_g \phi
(\tau,\cdot)\|_{L^\infty(D_\tau)} +\sum_{|I|\leq 2} (1+\tau)^{-1} \|
\, Z^I \phi(\tau,\cdot)\|_{L^\infty(D_\tau)}\Big)\, d\tau,
\end{multline}
If $H\!=\!0$, \eqref{eq:waveframeintro} is the decomposition in
radial and spherical coordinates and
\eqref{eq:strongdecaywaveeqintro} was used in \cite{L1}.

In section \ref{decayfirstder} we use the weak decay estimates \eqref{eq:weakdecayintro} in
\eqref{eq:strongdecaywaveeqintro} to get
the {\bf sharp decay estimates}
\beq\label{eq:strongdecay0intro}
|\pa \phi|\leq c_1\varepsilon (1+t)^{-1},\qquad\qquad
|\phi|\les c_1\varepsilon (1+|\, t-r|)(1+t)^{-1}.
\eq
The last inequality follows by integrating the first from $r\!=\!t\!+\!1$ where $\phi$ vanishes.
If \eqref{eq:strongdecay0intro}
had been true also for $Z\phi$ it would have been easy, but there is a small loss that requires a delicate analysis.

With the sharper decay \eqref{eq:strongdecay0intro} for $H(\phi)$
the decomposition of the wave operator \eqref{eq:waveframeintro} simplifies to
\beq \label{eq:waveframe2intro}
\Big|L_2^\alpha\pa_\alpha( p \,\pa_q \phi)-r\Boxr_g\phi\Big|
 \leq\frac{C}{1\!+t}\,\,\sum_{|I|\leq 2} \,|Z^I\phi|,
 \eq
where $p=r+t$ and
 \beq\label{eq:L2defintro}
L_2^\alpha\pa_\alpha=\big(L^\alpha-H^{\u L\u L} L^\alpha\big)\pa_\alpha=(\pa_t+\pa_r) +H^{\u L\u L}(\pa_r-\pa_t).
\eq

In section \ref{eikonal} we study the integral curves of the vector field
\eqref{eq:L2defintro} since we will integrate
\eqref{eq:waveframe2intro}.
Let
$q=r-t$, $p=r+t$ and $\omega=x/|x|$, and introduce the radial characteristics
$q=q(s,\rho,\omega)$ by
\beq\label{eq:qchardefintro}
d q/ ds =2H^{\u L\u L},
\qquad \text{when }\,\, |\,t-r|\leq t/2,\qquad q=\rho, \quad\text{when }\,\, |\, t-r|\geq t/2,\qquad p=2s.
\eq
Equivalently let $\rho$ be the solution of a radial eikonal equation:
 \beq\label{eq:rhodefintro}
L_2^\alpha\pa_\alpha \rho=0,\quad\text{when}\quad |\,t-r|\leqq
t/2,\qquad \rho=r-t, \quad\text{when} \quad|\,t-r|>t/2.
 \eq
$\rho$ behaves roughly like $q$:
\beq\label{eq:rhoestintro}
 \Big(\frac{1+t}{1\!+\!|\,\rho|}\Big)^{-c_1\varepsilon}
  \!\!\leq \frac{\pa \rho}{\pa q}
\leq  \Big(\frac{1+t}{1\!+\!|\,\rho|}\Big)^{c_1\varepsilon},
\qquad\quad
  \Big(\frac{1+t}{1\!+\!|\,\rho|}\Big)^{-c_1\varepsilon}\!\!\leq \frac{1+|\rho|}{1+|q|}
\leq  \Big(\frac{1+t}{1\!+\!|\,\rho|}\Big)^{c_1\varepsilon}.
\eq
This comes from differentiating \eqref{eq:rhodefintro}:
$L_2^\alpha\pa_\alpha \pa_q\rho+2\pa_q H^{\u L\u L}\, \pa_q \rho=0$,
and multiplying by the integrating factor using the estimate \eqref{eq:strongdecay0intro} for $H(\phi)$, as was observed in the
spherically symmetric case in \cite{L2}.

In section \ref{decaysecondder} we prove the following {\bf sharp decay estimates} for second derivatives:
\beq\label{eq:strongseconderestintro}
|\pa \phi| \leq \frac{c_1\varepsilon}{1+t}\frac{1}{(1+|\rho|)^\nu},\qquad\qquad
|\pa^2 \phi|\leq \frac{c_2\varepsilon}{1+t}\Big|\frac{\pa \rho}{\pa q}\Big|\frac{1}{(1+|\rho|)^{1+\nu}},\qquad \nu>0.
\eq
The first estimate follows from integrating \eqref{eq:waveframe2intro} along the integral curves of $L_2$
from $t\!=\!2|\rho|$ using
\eqref{eq:weakdecayintro}.
For the proof of the second we note that since $[\pa_\rho, L_2^\alpha\pa_\alpha]=0$ it follows from
\eqref{eq:waveframe2intro} (c.f. \cite{L2})
\beq\label{eq:seconderestrhointro}
\Big|L_2^\alpha\pa_\alpha( p \,\pa_\rho \pa_q \phi)-r\pa_\rho \Boxr_g\phi\Big|
 \leq\frac{C\rho_q^{-1}}{1+t+r}\sum_{|I|\leq 2} \,|\pa Z^I\phi|,
\qquad\text{where}\quad \pa_\rho =\rho_q^{-1}\pa_q.
\eq
The second estimate in \eqref{eq:strongseconderestintro}
follows from integrating
\eqref{eq:seconderestrhointro} using \eqref{eq:Zpaq} and \eqref{eq:rhoestintro}.

For vector fields we are not quite as lucky and there is a loss in the
{\bf strong decay estimate}:
\beq\label{eq:onevectorfieldestintro}
|Z\phi|\leq c_2 \varepsilon (1+t)^{-1+c_2\varepsilon}
(1+|q|)^{1-\nu}
\leq \frac{c_2\varepsilon}{1+t}\Big( (1+|q|) + (1+t)^{c_2\varepsilon/\nu}\Big).
\eq
In fact, by \eqref{eq:strongseconderestintro} the
commutator is
\beq
\sum_{|I|\leq 1} \big|\Boxr_g Z^I\phi\big|\leq C\sum_{|I|\leq 1} |Z^I H|\,|\pa^2 \phi|
\leq \frac{c_2^\prime\varepsilon \rho_q }{1+t}\,
\sum_{|I|\leq 1} \frac{|Z^I\phi|}{1+|\rho|}.
\eq
If we use \eqref{eq:waveframe2intro} applied to $Z^I\phi$ in place of $\phi$ and also commute $\rho_q^{-1}$ through \eqref{eq:waveframe2intro}, we get
an extra term due to that
 $L_2^\alpha\pa_\alpha \rho_q=-2\rho_q  \pa_q H^{\u L \u L}$:
\beq\label{eq:L2onevectorfield}
\sum_{|I|\leq 1} \Big|L_2^\alpha\pa_\alpha( p\, \pa_\rho Z^I\phi)\Big|
\leq
c_2\varepsilon \sum_{|I|\leq 1} \Big(\frac{|Z^I \phi|}{1+|\rho|} + |\pa_\rho Z^I\phi|\Big)
+ \frac{C\rho_q^{-1} }{1+t+r}\sum_{|I|\leq 3} \,|Z^I \phi|.
\eq
Since $\phi=0$ when $\rho\leq -1$ we can estimate $|Z^I\phi|/(1+|\rho|)$ by the derivative $|\pa_\rho Z^I\phi|$ we get if we first integrate
\eqref{eq:L2onevectorfield} from $t=2|\rho|$ where we can use
\eqref{eq:weakdecayintro} and \eqref{eq:Zpaq}:
\beq
M(t)\leq\int_1^t \frac{c_2\varepsilon}{1+\tau} M(\tau)\, d\tau +c_0\varepsilon
\qquad\text{where} \quad
M(t)=(1+t)\sup_{|\rho|\leq t/2} (1+|\rho|)^\nu \sum_{|I|\leq 1}
 |\pa_\rho Z^I\phi|
\eq
which by a Gronwall type of argument implies that
$M(t)\leq c_0 (1+t)^{c_2\varepsilon}$ from which \eqref{eq:onevectorfieldestintro} follows.

For more vector fields there is a problem with the most straightforward approach. We have
\beq
 |\Boxr_g Z^I\phi|\les |Z^I\phi|\, |\pa^2\phi|
+\sum_{|J|\leq 1,\,\, |K|=|I|-1\!\!\!\!\!\!\!\!\!\!\!\!\!\!
\!\!\!\!\!\!\!\!\!\!\!\!\!\!\!\!\!\!\!\!\!} |Z^J \phi|\, |\pa^2 Z^K\phi|
+\!\!\sum_{|J|+|K|\leq |I|,\,\, |J|<|I|,\,\, |K|<|I|-1\!\!\!\!\!\!\!\!\!\!\!\!\!\!\!\!\!\!\!\!\!\!\!\!\!\!\!\!\!\!\!\!\!\!\!\!
\!\!\!\!\!\!\!\!\!\!\!\!\!\!} |Z^J \phi|\, |\pa^2 Z^K\phi|.
 \eq
The first term can be handled as above and the terms in the last are
lower order. However the problem is the term with $|K|=|I|-1$ and
$|J|=1$ which is highest order. Using \eqref{eq:Zpaq} and
\eqref{eq:onevectorfieldestintro}
\beq\label{eq:highervectorfieldcommutatorintro} |Z\phi|
\sum_{|K|=|I|-1} |\pa^2 Z^K\phi| \leq \frac{c_2\varepsilon}{1+t}
\sum_{|J|=|I|} |\pa Z^J\phi|+
\frac{c_2\varepsilon}{(1+t)^{1-c_2\varepsilon/\nu}} \sum_{|K|=|I|-1}
|\pa^2 Z^K\phi|. \eq The first sum can be handled as above, however
the lack of decay in the last sum cause a problem. The estimate
\eqref{eq:onevectorfieldestintro} can not be improved to get the
needed $\varepsilon(1+t)^{-1}$ decay close to the light cone. One
could use modified vector fields that take into account the bending
of the light cones at infinity. The modified rotations are defined
by $(\widetilde{\Omega} \phi)(q,p,\omega)\!=\!\Omega
(\phi(q(\rho,s,\omega),2s,\omega))$, where $q(\rho,s,\omega)$ is as
in \eqref{eq:qchardefintro}. This  however leads to the loss of
regularity encountered by \cite{A2} for the energy estimate. We will
take a different approach. We can handle a $t^{c\varepsilon}$ loss
in terms of quantities we have already estimated. Therefore we will
first estimate $\pa^2 Z^{K} \phi$ in
\eqref{eq:highervectorfieldcommutatorintro} before we estimate $\pa
Z^I\phi$. In estimating $\pa^2 Z^K \phi$ a commutator term of the
same form  shows up with $\pa^2 Z^K\phi$ replaced by $\pa^3
Z^L\phi$, where $|L|\!=\!|K|\!-\!1$, so we must first estimate
$\pa^3 Z^L \phi$ and so on until finally we are left with
$\pa^{1+|I|} \phi$.

In section \ref{sharpdecayfour} we use induction to prove the following {\bf sharp decay estimate} for higher derivatives:
\beq\label{eq:stronghigherderestintro}
|\pa^{{k}} \phi|\leq
\frac{c_k\varepsilon}{1+t}
\Big(\frac{1+t}{1\!+\!|\,\rho|}\Big)^{c_k\varepsilon}\!\!
\frac{1}{(1+|\rho|)^{k-1+\nu}}.
\eq
In fact,  the commutator is
 \beq\label{eq:boxmanyder}
 |\Boxr_g \pa^{\,n} \phi|\leq  C |\pa \phi|\, |\pa^{1+n} \phi|
 +C\!\!\!\!\!\sum_{{k}_1+\cdots+{k}_\ell={n}+2,\,
 1\leq {k}_j\leq {n}, \, \ell\geq 2
 \!\!\!\!\!\!\!\!\!\!\!\!\!\!\!\!\!\!\!\! \!\!\!\!\!\!\!\!\!\!\!\!\!\!\!\!\!\!\!\! \!\!\!\!\!\!\!\!\!\!\!\!\!}
 |\pa^{\,k_1} \phi|\cdots |\pa^{{\,k}_{\,\ell}} \phi|
 \eq
where the first term has sufficient decay since $|\pa\phi|\leq c_1\varepsilon(1+t)^{-1}$ and the second sum is lower order.

Finally in section \ref{strongdecaymanyvectorfields} we use induction in $k$
as described above to show that
\beq\label{eq:weightedstrongdecaykIcopyintro}
|\pa^{i-k} Z^K \phi|\leq
c_{k,i}\,\varepsilon(1+t)^{-1+c_{k,i}\varepsilon}(1+|\rho|)^{1-(i-k)-\nu}\!\!\!,
\qquad i=|I|, \,\, k=|K].
 \eq

\section{Expressing the metric and the wave operator in the null frame}
\label{nullframe}
We introduce a nullframe for the Minkowski metric, ${\cal U}=\{
L,\u{L},S_1,S_2\}$, where
 \beq
 L^0=1,\quad L^i=\omega^i
 ,\qquad \u{L}^0=1,\quad
 \u{L}^i=-\omega^i,\qquad \omega^i=\frac{x^i}{|x|}, \qquad
 i=1,2,3,
 \eq
 and $S_1$ and $S_2$ are two smooth orhtonormal vector fields
 on the tangent space of
 the sphere $T(S^2)$. (We remark that these only exist locally so
 one has to work in a coordinate chart.)
 We will raise and lower the indices with respect to the Minkowski metric
 $V_\alpha=m_{\alpha\beta} V^\beta$, $V^\alpha=m^{\alpha\beta}
 V_\beta$, $m^{\alpha\beta}=m_{\alpha\beta}$.
(Here $m_{00}=-1$, $m_{ii}=1$, $i=1,2,3$, $m_{\alpha\beta}=0$, if $\alpha\neq \beta$.)
 We can express a vector field $X$ or the corresponding
 one form in the nullframe
 \beq\label{eq:decomp1}
 X^\alpha=X^L L^\alpha+X^{\u L} \u L^\alpha +X^A A^\alpha,\qquad
 X_\alpha=X^L L_\alpha+X^{\u L} \u L_\alpha +X^A A_\alpha.
 \eq
Here and in what follows $A,B,C,...$ denotes any of the vectors $S_1,
S_2$, and we used the convention that we sum over repeated upper and lower indices;
 \beq
X^A A^\alpha= X^{S_1}
S_1^\alpha+X^{S_2} S_2^\alpha .
 \eq
 The components can be calculated from the contractions:
\beq\label{eq:comp1}
X^L=-\frac{1}{2} X_{\underline{L}},\quad X^{\u L}=-\frac{1}{2} X_{L}, \quad
X^A=X_A,\qquad\text{where}\qquad X_Y=X_\alpha Y^\alpha=m_{\alpha\beta}X^\alpha Y^\beta=Y_X.
\eq
(This follows since $L^\alpha \underline{L}_\alpha=-2$,
$L^\alpha L_\alpha=\underline{L}^\alpha \underline{L}_\alpha=0$,
$L^\alpha A_\alpha=\underline{L}^\alpha A_\alpha=0$,
and $A^\alpha B_\alpha=\delta_{AB}$.)
  Recall that the inverse of the Minkowski metric $m^{\alpha\beta}$ can be expressed
in a nullframe
$$
m^{\alpha\beta}=-\frac{1}{2}\big(L^\alpha\underline{L}^\beta+
 \underline{L}^\alpha
 L^{\beta}\big)+\delta_{AB} A^\alpha B^\beta,
 $$
 where $\delta_{AB} A^\alpha B^\beta=S_1^\alpha S_1^\beta+S_2^\alpha
 S_2^\beta$.
 We make a similar decomposition for the bilinear form $g^{\alpha\beta}$:
 \beq\label{eq:decomp2}
 g^{\alpha\beta}=  g^{UV} U^\alpha V^\beta.
 \eq
 Here and in what follows $U,V,W,...$ denotes any vector in
 ${\cal U}=\{L,\u L, S_1,S_2\}$, and we used the convention that we sum over repeated upper and lower indices.
 The components can be calculated in terms of the
 contractions as follows:
 \beq\label{eq:comp2}
g^{L\u L}\!=\frac{1}{4} g_{\u L L},\quad g^{LL}\!=\frac{1}{4} g_{\u L \u
L},\quad g^{\u L \u L}\!=\frac{1}{4} g_{LL} ,\quad g^{LA}\!=-\frac{1}{2} g_{\u L A},\quad
 g^{\u L A}\!=-\frac{1}{2} g_{L A},\quad  g^{AB}\!\!=g_{AB},
\eq
where
\beq
g_{UV}=g^{\alpha\beta}U_\alpha V_\beta=g_{\alpha\beta} U^\alpha V^\beta,
\qquad\text{if} \quad  g_{\alpha\beta}=m_{\alpha\alpha^\prime} m_{\beta\beta^\prime} g^{\alpha^\prime\beta^\prime}
 \eq
denotes the lowering of indices with respect to the Minkowski metric
and not the inverse of $g^{\alpha\beta}$.
\begin{lemma} Suppose that $g^{\alpha\beta}$ is symmetric. Then
 \beq\label{eq:metricframe}
 g^{\alpha\beta}=-\frac{1}{2}\big(\, L_1^\alpha \underline{L}^\beta+
 \underline{L}^\alpha L_1^{\beta}\big)+\gamma^{\alpha\beta},
 \eq
 where
 \beq
 L_1^{\alpha}
 =g^{\mu\alpha} L_{\mu}+\frac{1}{4}g^{\mu\nu}L_\mu L_\nu
 \underline{L}^{\alpha}=-\frac{1}{2}g_{L\u L} L^\alpha
-\frac{1}{4} g_{LL}\u L^\alpha +g_{L}^{\,\, A} A^\alpha,
 \eq
 \beq
 \gamma^{\alpha\beta}=g^{LL}L^\alpha L^\beta+g^{AL} A^\alpha L^\beta+ g^{AL} L^\alpha A^\beta +g^{AB} A^\alpha B^\beta.
  \eq
 \end{lemma}
\begin{proof} Using  \eqref{eq:decomp1} and \eqref{eq:decomp2} we can write
\begin{align}
g^{\alpha\beta}&=g^{\u L \u L} \u L^\alpha \u L^\beta + g^{L\u L} L^\alpha \u L^\beta
+g^{ \u L L} \u L^\alpha L^\beta+ g^{\u L A} \u L^\alpha A^\beta +
g^{A \u L} A^\alpha \u L^\beta+\gamma^{\alpha\beta},\\
g^{\alpha\beta}&=g^{\u L \beta} \u L^\alpha + g^{\alpha\u L} \u L^\beta
-g^{ \u L \u L} \u L^\alpha \u L^\beta +\gamma^{\alpha\beta},
\end{align}
and the lemma follows from using  \eqref{eq:comp1} and  \eqref{eq:comp2}.

\end{proof}

 Let us introduce some further notation
 \beq
 \pa_q=\frac{1}{2}(\pa_r-\pa_t)=-\frac{1}{2}\u L^\alpha\pa_\alpha,
\qquad
 \pa_p=\frac{1}{2}(\pa_r+\pa_t)=\frac{1}{2} L^\alpha\pa_\alpha,
 \eq
 and
\beq
 \pab=\{\pa_L,\pa_{S_1},\pa_{S_2}\},\quad
 \text{where}\quad \pa_Y=Y^\alpha\pa_\alpha.
 \eq

\begin{lemma} Suppose that $g^{\alpha\beta}$ is a symmetric and bounded.
 Then
\beq\label{eq:derHLLderq}
\big| g^{\alpha\beta}\pa_\alpha\pa_\beta\phi -g_{LL}\pa_q^2 \phi \big|
\leq C|\overline{\pa} \pa\phi|
\eq
and with $\overline{\operatorname{tr}}\,
g=\delta_{AB} g^{AB}$ we have
 \beq \label{eq:waveframe}
\Big| 2 L_1^{\alpha}\pa_\alpha \pa_q \phi
+\frac{\overline{\operatorname{tr}}\, g}{r}\pa_q
\phi-g^{\alpha\beta}\pa_\alpha\pa_\beta \phi\Big|\leq C\sum_{1\leq |k|\leq
2\!\!\!\!\!\!\!\!} r^{|k|-2}\,|\pab^{\,k}\phi| .
 \eq
 Furthermore
 \beq\label{eq:eikonalframe}
\Big|g^{\alpha\beta}\pa_\alpha\rho\, \pa_\beta \rho +(\pa_q
\rho )L_1^{\alpha}\pa_\alpha \rho
-\delta^{AB}\pa_A\rho\,\, \pa_B\rho\Big|\les
\big(|H_{LL}|+|H_{LA}|+|H_{AB}|\big) |\overline{\pa}\rho|^2,
 \eq
where $H^{\alpha\beta}=g^{\alpha\beta}-m^{\alpha\beta}$.
\end{lemma}
\begin{proof} \eqref{eq:derHLLderq} and \eqref{eq:eikonalframe}
follow directly from \eqref{eq:metricframe}. By \eqref{eq:metricframe}
\beq
g^{\alpha\beta}\pa_\alpha\pa_\beta\phi=-L_1^\alpha \u
L^\beta\pa_\alpha\pa_\beta\phi
+\gamma^{\alpha\beta}\pa_\alpha\pa_\beta\phi.
\eq
Now
\beq
-L_1^\alpha \u L^\beta\pa_\alpha\pa_\beta\phi =-L_1^\alpha
\pa_\alpha\big(\u L^\beta\pa_\beta\phi\big) +\big(L_1^\alpha
\pa_\alpha \u L^\beta\big)\pa_\beta\phi .
\eq
Note that if
$\omega_j=x_j/|x|$ then \beq
\pa_i\,\omega_j=\frac{1}{r}\big(\delta_{ij}-\omega_i \,
\omega_j\big), \qquad\text{so}\qquad L^i\pa_i\,\omega^j=\u L^i\pa_i
\,\omega^j=0,\qquad A^i\pa_i\,\omega^j=\frac{1}{r} A^j. \eq Hence $
L_1^\alpha\pa_\alpha \u L^\beta=-L_1^A A^\beta /r=-g_{L}^{\,\, A} A^\beta/r$ and it follows that
\beq
-L_1^\alpha \u L^\beta\pa_\alpha\pa_\beta\phi =2
L_1^\alpha\pa_\alpha \pa_q+g_{L}^{\,\,\, A} r^{-1} A^\alpha
\pa_\alpha\phi.
\eq
We have
\begin{multline}
\gamma^{\alpha\beta}\pa_\alpha\pa_\beta\phi=g^{LL}L^\alpha L^\beta\pa_\alpha\pa_\beta\phi+2g^{LA} A^\alpha L^\beta\pa_\alpha\pa_\beta\phi+g^{AB} A^\alpha B^\beta\pa_\alpha\pa_\beta\phi\\
=g^{LL}L^\alpha \pa_\alpha(L^\beta\pa_\beta\phi)+2g^{LA} A^\alpha \pa_\alpha(L^\beta\pa_\beta\phi)+g^{AB} A^\alpha \pa_\alpha
(B^\beta\pa_\beta\phi)-g^{AB} (A^\alpha\pa_\alpha B^\beta)\pa_\beta\phi
\end{multline}
Since $\omega_j B^j=0$ it follows that \
\beq
(A^i\pa_i B^j)\omega_j=-A^i B^j\pa_i\omega_j=-A_j B^j\frac{1}{r}=-\frac{1}{r}\delta_{AB},\qquad
\pa_i=\omega_i\pa_r+\pab_i,
\eq
where $\pab_i$ is tangential
it follows that
\beq
-g^{AB} (A^\alpha\pa_\alpha B^\beta)\pa_\beta\phi=g^{AB}\delta_{AB}\,
\frac{1}{r}\pa_r \phi-g^{AB} (A^i\pa_i B^j)\pab_j\phi.
\eq
Since $B^j$ are smooth functions of $\omega$ it follows that
$|A^i\pa_i B^j|\leq C/r$. The lemma therefore follows from the above identities.
\end{proof}

\section{The vector fields associated with the wave operator, commutators.}\label{vectorfield}
Let $Z\in {\cal Z}$ be any of the vector fields
$$
\Omega_{\a\b} = -x_{\a}\pa_{\b} + x_{\b}\pa_{\a},\qquad
 S = t\pa_{t} + r\pa_{r},\qquad \pa_\alpha,
$$
where $x_0=-t$ and $x_i=x^i$, for $i\geq 1$. Let
$I=(\iota_1,...,\iota_k)$, where $|\iota_i|=1$, be an ordered
multiindex of length $|I|=k$ and let $Z^I=Z^{\iota_1}\cdot\cdot\cdot
Z^{\iota_k}$ denote a product of $|I|$ such derivatives. With a
slight abuse of notation we will also identify the index set with
vector fields, so $I=Z$ means the index $I$ corresponding to the
vector field $Z$. Furthermore, by a sum over $I_1+I_2=I$ we mean a
sum over all possible order preserving partitions of the ordered
multiindex $I$ into two ordered multiindices $I_1$ and $I_2$, i.e.
if $I=(\iota_1,...,\iota_k)$, then
$I_1=(\iota_{i_1},...,\iota_{i_n})$ and
$I_2=(\iota_{i_{n+1}},...,\iota_{i_k})$, where $i_1,...,i_k$ is any
reordering of the integers $1,...,k$ such that $i_1<...<i_n$ and
$i_{n+1}<...<i_k$ and $i_1,...,i_k$.
 With this convention
Leibnitz rule becomes $Z^I(fg)=\sum_{I_1+I_2=I} (Z^{I_1} f)(Z^{I_2}
g)$.

We recall that the family ${\cal Z}$ possesses  special commutation
properties: for any vector field $Z\in {\cal Z}$
$
[Z,\Box]=-C_Z \Box,
$
where the constant $C_Z$ is only different from zero in the case of
the scaling vector field $C_S=2$. Moreover $[Z,\pa_\alpha]=C_{Z\alpha}^{\,\beta}\pa_\beta$,
for some constants $C_{Z\alpha}^{\,\beta}$.
It is easy to show the following
identities
\beq\label{eq:iZ}
\pa_{t} = \frac {t S -
x^{i}\Omega_{0i}}{t^{2}-r^{2}}, \qquad \pa_{i} = \frac
{-x^{j}\Omega_{ij} + t \Omega_{0i} - x_{i} S}{t^{2}-r^{2}},
 \eq
 and for some smooth functions $f_A^{ij}(\omega)$;
 \beq \label{eq:somega}
 \pa_{L} =
\frac {S + \omega^{i}\Omega_{0i}}{t+r},\qquad
\pa_A=\frac{f_A^{ij}(\omega)\,\Omega_{ij}}{r},\qquad
\Omega_{ij}=\frac{r}{t}\big(\omega_j\Omega_{0i}-\omega_i\Omega_{0j}\big).
\eq
Recall that $\bar\pa$ denotes the tangential derivatives, i.e.,
$\pa_T$, where $T\in\, {\cal T}=\{L, {S_1},{S_2}\}$.
\begin{lemma}\label{tanderZ}
For any function $\phi$;
\begin{align}
&(1+t+r)|\bar \pa  \phi|+(1+|t-r|)|\pa \phi|\leq C
\sum_{|I|=1}|Z^I \phi|,\label{eq:tanZ}\\
&(1+t+r)|\pa  \phi|\leq C\, r|\pa_q \phi|+C
\sum_{|I|=1}|Z^I \phi|,\label{eq:derZ}\\
&|\bar\pa^2 \phi|+r^{-1}|\pab \phi|\leq \frac{C}{r}\sum_{|I|\leq 2}
\frac{|Z^I \phi|}{1+t+r}, \quad{\text {where}} \quad |\bar\pa^2
\phi|^2=\sum_{S,\, T\in\, {\cal T}} |\pa_S \pa_T \phi|^2,
\label{eq:2tanZ}\\
&(1+|t-r|)^{k}|\pa^k \phi|\leq C
\sum_{|I|\leq |\bold{k}|}|Z^I \phi|.\label{eq:derframeZ}
\end{align}
\end{lemma}
\begin{proof} First we note that if $r+t\leq 1$ then
\eqref{eq:tanZ} holds since the standard derivatives $\pa_\alpha$
are included in the sum on the right. The inequality for $|\bar\pa
\phi|$ in \eqref{eq:tanZ} follows directly from \eqref{eq:somega}.
The inequality for $|\pa \phi|$ in \eqref{eq:tanZ} follows from \eqref{eq:iZ}.
The inequality \eqref{eq:derZ} follows similarly from \eqref{eq:iZ}-\eqref{eq:somega}.
The proof of \eqref{eq:2tanZ} follows immediately from
\eqref{eq:somega} and the inequality $|\pa_i\omega_j|\leq C r^{-1}$.
The inequality \eqref{eq:derframeZ} follows from
repeated use of \eqref{eq:tanZ} and the commutator identity
$[Z,\pa_i]=c_i^\alpha\pa_\a$, where $c_i^\a$ are constants.
\end{proof}

Let $H^{\alpha\beta}=g^{\alpha\beta}-m^{\alpha\beta}$ and
 \begin{align}\label{eq:intvf}
L_1^\alpha&=L^\alpha-\frac{1}{2}H_{L\u L} L^\alpha -\frac{1}{4} H_{LL}\u
L^\alpha +H_{L}^{\,\, A} A^\alpha, \\
\label{eq:intvf2}
L_2^\alpha\pa_\alpha &=(L^\alpha-\frac{1}{4}H_{LL}\u L^\alpha)\pa_\alpha=2\pa_p+\frac{1}{2}H_{LL} \pa_q.
 \end{align}

\begin{lemma}\label{wavei} Let $\Boxr_g=g^{\alpha\beta}\pa_\alpha\pa_\beta$,
$H^{\alpha\beta}=g^{\alpha\beta}-m^{\alpha\beta}$ and suppose that
$|H|\leq 1/4$. Then
 \beq \label{eq:phiwave1}
\Big|\Big( 2
L_1^{\alpha}\pa_\alpha+\frac{\ell}{r} \Big)( r \pa_q
\phi)-r\Boxr_g\phi\Big|\leq \frac{C}{1+t+r}\!\sum_{|I|\leq 2}
|Z^I\phi|, \qquad\quad \ell=\overline{\operatorname{tr}}\, H
\!+ H_{L\Lb}\!-\tfrac{1}{2}H_{LL},
\eq
where $\overline{\operatorname{tr}}\, H=\delta^{AB} H_{AB}$.
 Suppose also that
\beq\label{eq:wavecoorddecay}
 |H_{LL}|+|H_{LA}|+|H_{AA}|+|H_{L\underline{L}}|\leq
\frac{1+|\,t-r|}{1+t+r}\,\,\Big|\frac{1+t+r}{1+|\,t-r|}\Big|^a,\qquad a\geq 0.
 \eq
 Then \beq \label{eq:phiwave2}
\Big|2L_2^\alpha\pa_\alpha( r \pa_q \phi)-r\Boxr_g\phi\Big|
 \leq\frac{C}{1+t+r}\,\,\Big|\frac{1+t+r}{1+|\,t-r|}\Big|^a\sum_{|I|\leq 2} \,|Z^I\phi|.\eq
\end{lemma}
\begin{proof} \eqref{eq:phiwave1} follows from \eqref{eq:waveframe} using
\eqref{eq:2tanZ} and
\beq
\overline{\operatorname{tr}}\, g=2+\overline{\operatorname{tr}}\, H,\qquad
2L_1^\alpha\pa_\alpha r=2-H_{L\u L}+\frac{1}{2} H_{LL}.
\eq
\eqref{eq:phiwave2} follows from \eqref{eq:phiwave1} and \eqref{eq:tanZ} using that
\beq
r |L^\alpha\pa_\alpha (\u L^\beta\pa_\beta \phi)|= r |L^\alpha \u L^\beta \pa_\alpha\pa_\beta\phi|\les \sum_{|I|=1, \,\beta=0,1,2,3\!\!\!\!\!\!\!\!\!\!\!\!\!\!\!\!\!\!\!\!} |Z^I \pa_\beta\phi|
\les \sum_{|I|\leq 1} |\pa Z^I \phi|\les \frac{1}{1+|\, t-r|} \sum_{|I|\leq 2\!\!\!}
|Z^I \phi|.
\eq
\end{proof}

\begin{lemma}\label{eq:secondderest}
Suppose that $H$ satisfy the assumptions of Lemma \ref{wavei} and
\beq\label{eq:wavecoorddecay}
 |\pa H|\leq
\frac{1}{1+t+r}\,\,\Big|\frac{1+t+r}{1+|\,t-r|}\Big|^a,\qquad a\geq 0.
 \eq
Then
\beq\label{eq:secondereq}
\Big|2L_2^\alpha\pa_\alpha( r \pa_q^2 \phi)+r (\pa_q H_{LL})
\pa_q^2\phi -r\pa_q \Boxr_g\phi\Big|
 \leq\frac{C}{1+t+r}\,\,\Big|\frac{1+t+r}{1+|\,t-r|}\Big|^a\sum_{|I|\leq 2} \,|\pa Z^I\phi|.
\eq
\end{lemma}
\begin{proof}To prove \eqref{eq:secondereq} we first commute
$\pa_\gamma \Boxr_g \phi=\Boxr_g \pa_\gamma\phi+(\pa_\gamma H^{\alpha\beta})\pa_\alpha\pa_\beta\phi$,
use \eqref{eq:derHLLderq} applied to $\pa_\gamma\, H^{\alpha\beta}$;
\beq
\big| (\pa_\gamma H^{\alpha\beta})\pa_\alpha\pa_\beta\phi-(\pa_\gamma H_{LL})\pa_q^2\phi\big|
\les |\pa H|\, |\pab\pa\phi|\les |\pa H|\, \sum_{|I|=1}|\pa Z^I \phi|.
\eq
Using \eqref{eq:phiwave2} applied to $\pa_\gamma\phi$ in place of $\phi$ now gives
\beq\
\Big|2L_2^\alpha\pa_\alpha( r \pa_q \pa_\gamma\phi)+r (\pa_\gamma H_{LL})
\pa_q^2\phi -r\pa_\gamma \Boxr_g\phi\Big|
 \leq\frac{C}{1+t+r}\,\,\Big|\frac{1+t+r}{1+|\,t-r|}\Big|^a\sum_{|I|\leq 2} \,|\pa Z^I\phi|.
\eq
and contracting with $\u L^\gamma$, using that it commutes with $L_2$ and $\pa_q$, gives \eqref{eq:secondereq}
\end{proof}

Let us now calculate the commutators of vector fields with $\Boxr_g=g^{\alpha\beta}\pa_\alpha\pa_\beta=\Box+H^{\alpha\beta}\pa_\alpha\pa_\beta$:
\begin{multline}
Z\Boxr_g \psi=Z\Box\psi+ Z\big( H^{\alpha\beta} \pa_\alpha\pa_\beta \psi\big)\\
=\Box Z\psi-C_Z \Box\psi+H^{\alpha\beta}\pa_\alpha\pa_\beta Z\psi+ 2 H^{\alpha\beta} C_{Z\alpha}^{\,\gamma} \pa_\gamma\pa_\beta\psi + (ZH^{\alpha\beta})\pa_\alpha\pa_\beta\psi\\
=\Boxr_g Z\psi -C_Z \Boxr_g\psi+ C_Z H^{\alpha\beta}\pa_\alpha\pa_\beta \psi
+2 H^{\alpha\beta} C_{Z\alpha}^{\,\gamma} \pa_\gamma\pa_\beta\psi + (ZH^{\alpha\beta})\pa_\alpha\pa_\beta\psi
\end{multline}
and hence with $\widehat{Z}=Z+C_Z$;
\beq\label{eq:onecommute}
\Boxr_g Z\psi=\widehat{Z}\Boxr_g \psi-
 \big( Z H^{\alpha\beta}+2 C_{Z\gamma}^{\,\alpha}  H^{\gamma\beta}+C_Z H^{\alpha\beta}\big)\pa_\alpha\pa_\beta\psi.
\eq
In general we have
\beq\label{eq:manycommute}
\Boxr_g Z^I\psi
=\widehat{Z}^I\Boxr_g\psi+\sum_{|J|+|K|\leq |I|,\,\,  |K|<|I|
\!\!\!\!\!\!\!\!\!\!\!\!\!\!\!\!\!\!\!\!\!\!\!\!\!\!}
C_{\!JK\, \alpha\beta}^{\,\, I\,\,\,\,\delta\gamma}\,(Z^{J} H^{\alpha\beta})\,\, \pa_\gamma\pa_\delta Z^{K}  \psi,
\eq
where $C_{\!JK\, \alpha\beta}^{\,\, I\,\,\,\,\delta\gamma}$ are constants.
The same formula holds for usual derivatives $\pa_\alpha$ in place of $Z$
even without the lower order terms with $|J|+|K|<|I|$, but we will need to separate these from the vector fields since they will behave better.
Let $\bold{k}=(k_1,...k_n)$ be a multindex and
$\pa^\bold{k}=\pa_{k_1}\cdots\pa_{k_n}$. Then
\beq\label{eq:manycommutemixed}
\Boxr_g\pa^{\bold{k}} Z^I\psi
=\pa^{\bold{k}} \widehat{Z}^I \Boxr_g\psi+\sum_{|J|+|K|\leq |I|,\, \bold{m}+\bold{n}=\bold{k},\,\, |\bold{n}|+|K|<|\bold{k}|+|I|
\!\!\!\!\!\!\!\!\!\!\!\!\!\!\!\!\!\!\!\!\!\!\!\!\!\!
\!\!\!\!\!\!\!\!\!\!\!\!\!\!\!\!\!\!\!\!\!\!\!\!\!\!
\!\!\!\!\!\!\!\!\!\!\!\!\!\!\!\!\!\!\!\!\!\!\!\!\!\!
\!\!\!\!\!\!\!\!\!\!\!}
 C_{\!JK\, \alpha\beta}^{\,\, I\,\,\,\,\delta\gamma}\,\,\,
(\pa^\bold{m} Z^J H^{\alpha\beta})\,\,\pa_\alpha\pa_\beta\pa^{\bold{n}} Z^K\psi .
\eq
Moreover
\beq\label{eq:Hder}
|\pa^{\bold{m}} Z^J H^{\alpha\beta}(\phi)|\leq
C\!\sum_{\bold{m}_1+\cdots+\bold{m}_\ell=\bold{m},\,\,
J_1+\cdots+J_\ell=J,\,\, \ell\geq 1
\!\!\!\!\!\!\!\!\!\!\!\!\!\!\!\!\!\!\!\!\!\!\!\!\!\!
\!\!\!\!\!\!\!\!\!\!\!\!\!\!\!\!\!\!\!\!\!\!\!\!\!\!
\!\!\!\!\!\!\!\!\!\!\!\!\!\!\!\!\!\!\!\!\!\!\!\!\!\!}
|\pa^{\bold{m}_1} Z^{J_1} \phi|\cdots
 |\pa^{\bold{m}_\ell} Z^{J_\ell} \phi|.
\eq
We have the following:
\begin{lemma} Suppose that $\Boxr_g\phi=0$ and $|\pa^\bold{n} Z^K\phi|\leq 1$,
for $|\bold{n}|+|K|\leq N-5$. Then for $|\bold{k}|+|I|\leq N$ we have
\beq\label{eq:bigcommuteest}
 |\Boxr_g \pa^{\bold{k}} Z^I\phi|\les
 \sum_{|\bold{n}|\leq |\bold{k}|,\, \, |J|+|K|\leq |I|,\,\,
|K|<|I|
 \!\!\!\!\!\!\!\!\!\!\!\!\!\!\!\! \!\!\!\!\!\!\!\!\!\!\!\!\!\!\!\!\!
 \!\!\!\!\!\!\!\!\!\!\!\!\!\!\!\!\!\!\!\!\!\!\!\!\!\!}   |Z^J \phi|\,\, |\pa^2 \pa^{\bold{n}} Z^K\phi|\,\,
+\sum_{|{\bold{m}}|+|\bold{n}|\leq |\bold{k}|,\,\,
 |J|+|K|\leq |I|
\!\!\!\!\!\!\!\!\!\!\!\!\!\!\!\!\!\!\!\!\!\!\!\!\!\!\!\!\!\!\!\!\!\!\!\!\!\!\!\!
\!\!\!\!\!\!\!\!\!\!\!\!\!}  |\pa \pa^{\bold{m}} Z^{J}\phi|\, |\pa \pa^{\bold{n}} Z^{K}\phi|.
\eq
If $|\bold{k}|=0$ then only the first sum is present and
if $|I|=0$ then only the second sum is present.
\end{lemma}
\begin{proof} If $|K|=|I|$ in the sum \eqref{eq:manycommutemixed} then
$|J|=0$, $|\bold{m}|\geq 1$ and $|\bold{n}|<|\bold{k}|$ so using \eqref{eq:Hder} we see that this term can be bounded by a term of the form in the second sum.
On the other hand if $|\bold{n}|=|\bold{k}|$ in \eqref{eq:manycommutemixed} then
$|\bold{m}|=0$ and $|K|<|I|$, and this term can be bounded by a term in the first sum above. Finally a term with $|\bold{n}|<|\bold{k}|$ and $|K|<|I|$ in \eqref{eq:manycommutemixed} can be bounded by a term contained in one of the sums above since under the assumptions of the lemma
\beq\label{eq:Hder2}
|\pa^{\bold{m}} Z^J H^{\alpha\beta}(\phi)|\leq
C\sum_{|\bold{n}|\leq |\bold{m}|,\,\,|K|\leq |J|
\!\!\!\!\!\!\!\!\!\!\!\!\!\!\!\!\!\!\!\!\!\!\!\!\!\!\!\!\!\!}
|\pa^{\bold{n}} Z^{K} \phi|\,.
\eq
\end{proof}

\section{Decay estimates for the wave equation on a curved background}
\label{strongdecaywave}
 For $(t,x)\in D=\{(t,x)\in\bold{R}\times\bold{R}^3;\, t/2<
|x|<3t/2\}$, let $X_i(s)\in\bold{R}^{1+3}$, $i=1,2$, be the backward integral curve
\beq\label{eq:intcurves}
\frac{d}{ds} X_i^\alpha=L_i^\alpha(X_i),\quad s\leq 0, \qquad X_i(0)=(t,x)
\eq
of the vector fields \eqref{eq:intvf}.
Let $s_i<0$ be the largest number such that $X(s_i)\in\pa D$, $X(s)\in D$, $s>s_i$. Let
$\tau_i=\tau_i(t,x)=X^0(s_i)$. Assuming that $|H|\leq 1/4$ the integral curve
will in fact intersect $\pa D$.

The following lemma is a generalization of a lemma in
\cite{L1,L-R1,L-R2}.

\begin{lemma} \label{strongdecaywaveeq} Suppose that $H^{\alpha\beta}=g^{\alpha\beta}-m^{\alpha\beta}$
satisfies $|H|\leq 1/16$ and either of the following
\beq\label{eq:decaymetric1} (1)\qquad\qquad\quad\qquad \int_0^T
\|H(t,\cdot)\|_{L^\infty(D_t)}\, \frac{dt}{1+t}\leq 1,\qquad\qquad
\qquad\qquad \eq
 where $D_t=\{x\in\bold{R}^3;\, t/2< |x|< 3t/2\}$, or
 \beq\label{eq:decaymetric2}
(2)\qquad\qquad  |H_{LL}|+|H_{LA}|+|H_{AA}|+|H_{L\underline{L}}|\leq
\frac{1}{4}\,\frac{1+|\, t-r|}{1+t+r}\,\,\Big|\frac{1+t+r}{1+|\,t-r|}\Big|^a, \qquad\text{in}\quad D,
 \eq
 where $D=\{(t,x)\in\bold{R}\times\bold{R}^3;\, t/2<|x|< 3t/2\}$
and $a\geq 0$.
 Then for any $a\geq 0$;
\begin{multline}\label{eq:strongdecaywaveeq}
(1+t+r) |\pa \phi(t,x)| \leq  C\sup_{\tau_i\leq \tau\leq \,t}
\sum_{|I|\leq 1}
\|Z^I\! \phi(\tau,\cdot)\|_{L^\infty}\\
+ C\int_{\tau_i}^t \Big((1+\tau)\| \Boxr_g \phi
(\tau,\cdot)\|_{L^\infty(D_\tau)} +\sum_{|I|\leq 2} (1+\tau)^{-1+a} \|
(1+|q(\tau,\cdot)|)^{-a}\, Z^I \phi(\tau,\cdot)\|_{L^\infty(D_\tau)}\Big)\, d\tau,
\end{multline}
where $q(t,x)=|x|-t$.
Here $\tau_i\!=\!\tau_i(t,x)$, where $i\!=\!1$ if condition (1) hold and $i\!=\!2$ if condition (2) hold, is defined as follows. If $(t,x)\notin D$
then $\tau_i=t$. If $(t,x)\in D$ then $\tau_i$ is the first time the
backward integral curve from $(t,x)$ for the vector field $L_i$,
in \eqref{eq:intvf}- \eqref{eq:intvf2},
leaves the region $D$. In general $0\leq \tau_i\leq t$.
\end{lemma}
\begin{proof}
By \eqref{eq:derZ} we only need to show that $\psi=r\, \pa_q \phi$ is
bounded by the right hand side.
Lemma \ref{wavei} can be summarized
 \beq
\Big|\Big( 2
L_i^{\alpha}\pa_\alpha\!+\!\frac{\ell_i}{r} \Big)\psi-r\Boxr_g\phi\Big|\leq \frac{C}{1\!+t\!+r}\Big|\frac{1+t+r}{1\!+|\,r\!-\!t|}\Big|^a\!\sum_{|I|\leq 2}
|Z^I\phi|,\quad\text{where}\quad
\begin{cases}  \ell_1\!\!\!\!\!\!&=\overline{\operatorname{tr}}\, H
\!+\!H_{L\Lb}\!-\!\tfrac{1}{2}H_{LL},\\
\ell_2\!\!\!\!\!\!&=0\end{cases}\!\!\!\!
\eq
With the integrating factor
\beq
G_i(s)=-2^{-1}\int_s^0 \ell_i\,(X_i(\sigma))/r(X_i(\sigma))\, d\sigma
\eq
we have along the integral curves \eqref{eq:intcurves}
\beq
\frac{d}{ds} \Big( \psi(X_i(s))e^{G_i(s)}\Big) =
\frac{1}{2} \, e^{G_i(s)}\Big( \big(2 L_i^\alpha\pa_\alpha +\frac{\ell_i}{r}\big)\psi\Big)
\big(X_i(s))
\eq
It follows from the assumption \eqref{eq:decaymetric1} that $|G_1|\leq C$ independently of $s$ and $G_2=0$. Hence it follows from integrating this from $s_i$ to $0$ that
 \begin{multline}
 |\psi(t,x)|\les |\psi(X_i(s_i))|+\int_{s_2}^0\Big|\Big(\big( 2 L_i^\alpha\pa_\alpha +\frac{\ell_i}{r}\big)\psi \Big)(X_i(s))\Big|\, ds\\
 \les \sum_{|I|\leq 1} |Z^I\phi(X_i(s_i))|+\int_{s_i}^0
 r|\Boxr_g\phi(X_i(s))|+\sum_{|I|\leq 2}(1+t+r)^{-1}|Z^I\phi(X_i(s))|\,
 ds.
 \end{multline}
 Since $t=X_i^0$ and $d X_i^0/ds= L_i^0$, where $L_2^0=1-\frac{1}{4}H_{LL}$ and
$L_1^0=1-\frac{1}{4}H_{LL}-\frac{1}{2}H_{L\u L}$ it follows that
$ 1/2\leq |\,d t/ds|\leq 2$ and the lemma follows.
\end{proof}

Next we define substitutes $\rho_i$ for $r-t$. Let
$\rho_i=\rho_i(t,x)$ be constant along the integral curves of $L_i$
and equal to $r-t$ outside a neighborhood of the forward light cone:
 \beq\label{eq:rhodefi}
L_i^\alpha\pa_\alpha \rho_i=0,\quad\text{when}\quad |\,t-r|\leqq
t/2,\qquad \rho_i=r-t, \quad\text{when} \quad|\,t-r|>t/2
 \eq
Note that $(t_i,x_i)$ is the first point the backward integral curve intersects
$|r-t|=t/2$ then $|\rho_i(t,x)|=|\,t_i-|x_i|\,|=t_i/2=\tau_i(t,x)/2\leq t/2$, since $t$ is increasing
along the forward integral curves.
Here $\tau_i$ was defined to be the smallest $t$ along the integral curve with
$|r-t|\leq t/2$. Hence
\beq
|\rho_i(t,x)|=\tau_i(t,x)/2\leq t/2,\qquad\text{when}\quad|\, t-|x|\, |\leq t/2
\eq
 We have
\begin{lemma}\label{weightedstrongdecaywaveeq}
 Suppose that either of the conditions in Lemma
\ref{strongdecaywaveeq} hold.
 Then for $\nu>\mu\geq 0$ and any $ a\geq 0$ we have
\begin{multline}\label{eq:weightedstrongdecaywaveeq} (1+t+|x|) (1+|\rho_i(t,x)|)^\nu |\pa
\phi(t,x)| \leq
 C\sup_{\tau_i\leq \tau\leq \,t} (1+\tau)^{\nu-\mu+b}
\sum_{|I|\leq 2} \| (1+|\rho_i(\tau,\cdot)|)^{\mu}
(1+|q(\tau,\cdot)|)^{-a}
Z^I\! \phi(\tau,\cdot)\|_{L^\infty}\\
+C\int_{\tau_i}^t (1+\tau)\| (1+|\rho_i(\tau,\cdot)|)^\nu\Boxr_g
\phi (\tau,\cdot)\|_{L^\infty(D_\tau)} \, d\tau,
\end{multline}
where $q(t,x)=|x|-t$.
Here $i=1$ if condition (1) holds and $i=2$ if condition (2) holds.
\end{lemma}
\begin{proof} This follows from Lemma \ref{strongdecaywaveeq} using that
$\rho_i(t,x)=\tau_i(t,x)/2\leq \tau/2$ along the integral curves and
$\rho_i$ are constant along the integral curves.
We also use that $\int_{2\rho_i}^\infty (1+\tau)^{-1-\nu+\mu}\, d\tau
=C (1+\rho_i)^{\mu-\nu}$.
\end{proof}

\section{Estimates for the radial characteristics and eikonal equation}
\label{eikonal}
We will use a curved substitute $\rho(t,x)$,
for the distance to the forward light cone $r-t$.
 Let
$\rho=\rho(t,x)$ be equal to $r-t$ outside a neighborhood of the forward light cone constant along the integral curves of the radial vector field $L_2$
 close to the light cone:
 \beq\label{eq:rhodef}
L_2^\alpha\pa_\alpha \rho=0,\quad\text{when}\quad |\,t\!-r|\leqq
t/2,\qquad \rho=r-t, \quad\text{when} \quad|\,t\!-r|\geq t/2,
 \eq
 where
 \beq
L_2^\alpha\pa_\alpha=2\pa_p +\frac{1}{2}H_{LL}\pa_q,
 \eq
and we think of $\rho=\rho(q,p,\omega)$ as a function of $q=r-t$, $p=r+t$ and $\omega=x/|x|$.
We call \eqref{eq:rhodef} the radial eikonal equation.
Alternatively, let $X_2(s)$ be the
integral curves of the vector field $L_2$, i.e.
$\dot{X}_2=L_2$.
 Then we can choose the initial conditions when $|\, t-r|=t/2$ so that
 \beq\label{eq:changeofvariables}
q=-X_{2 L}(s,\rho,\omega),\qquad\text{and}\qquad p=X_{2\u L}=2s,
\eq where
\beq\label{eq:qchar}
\frac{d}{d s} X_{2\u L}=-\frac{1}{2} H_{LL},
\qquad \text{when }\,\, |\,t-r|\leq t/2,\qquad q=\rho, \quad\text{when }\,\, |t-r|\geq t/2.
\eq
We call these the radial characteristics.

We now state the main estimate for $\rho$ assuming some estimates for $H_{LL}$ that will be proven later.
We will assume that $H_{LL}=0$, when $r>t+1$ and $t>$ so in fact
$\rho=r-t$, when $r>t+1$.

Since, as we show below,
$0<\pa \rho/\pa q<\infty$, $\rho$ is an invertible function of $q$ for fixed  $(p,\omega)$, satisfying $d\rho/dq =\pa \rho/\pa q$, and $q$ is an invertible function of $\rho$, satisfying $d q/ d\rho=(d \rho/dq)^{-1}$.
We have thus introduced a change of variables
$(\rho,s,\omega)\to(q(\rho,s,\omega),2s,\omega)$.
Note that multiplication by any function of $\rho$ commutes with $L_2^\alpha\pa_\alpha$ and a calculation using that shows that
\beq\label{eq:L2commuteparho}
[L_2^\alpha\pa_\alpha,\pa_q]=-\frac{\pa_q H_{LL}}{2}\pa_q,\quad\text{and}\quad [L_2^\alpha\pa_\alpha,\pa_\rho]=0,\qquad \text{if}\quad
\pa_\rho=\rho_q^{-1}\pa_q.
\eq
The following lemma was essentially proven in \cite{L1} in the spherically symmetric case:
\begin{prop}\label{eikonalone}
Let $\rho(t,x)$ be as in \eqref{eq:rhodef}
and suppose that $H_{LL}$ satisfies
\beq\label{eq:HLLassumption1}
|\pa H_{LL}|\leq \frac{c_1\varepsilon }{1\!+t}\,\frac{1}{(1\!+\!|\,\rho|)^{\nu}},
\quad\text{and}\qquad |H_{LL}|\leq c_1\varepsilon \, \frac{1\!+\!|\,q|}{1+t}
\eq
for some $\nu\geq 0$.
Then \beq \label{eq:eikonal1}
\Big(\frac{1+t}{1\!+\!|\,\rho|}\Big)^{-c_1\varepsilon V(\rho)}\!\!
\leq\frac{\pa \rho}{\pa q}
\leq\Big(\frac{1+t}{1\!+\!|\,\rho|}\Big)^{c_1\varepsilon V(\rho)} ,\qquad\quad V(\rho)=(1\!+\! |\,\rho|)^{-\nu}
 \eq
and
\beq \label{eq:eikonal2}
\Big(\frac{1+t}{1\!+\!|\,\rho|}\Big)^{-c_1\varepsilon}\!\! \leq \frac{1+|\,q|}{1+|\,\rho|}
\leq \Big(\frac{1+t}{1\!+\!|\,\rho|}\Big)^{c_1\varepsilon} \eq
\end{prop}
\begin{proof}
We have
\beq\label{eq:rhoqeq}
L_2^\alpha\pa_\alpha \pa_q\rho+\frac{1}{2}\pa_q H_{LL}\, \, \pa_q \rho=0.
\eq
Let $X(s)$ be a backward integral curve of the vector fields $L_2$:
\beq
\frac{d}{ds} X^\alpha=L_2^\alpha(X),\quad s\leq 0, \qquad X(0)=(t,x)
\eq
and let $s_2<0$ be the largest number such that
$X(s_2)=(t_2,x_2)$ satisfies $|t_2-|x_2|\,|=t_2/2$. If we multiply by the
integrating factor
\beq
G(s)=-\int_s^0 \frac{1}{2} (\pa_q H_{LL})(X(\tau))\, d \tau
\eq
we get
\beq
\frac{d }{ds } \Big(\frac{\pa \rho}{\pa q}\big( X(s)\big)e^{G(s)}\Big)=0.
\eq
Integrating this from $s_2$ (where $\pa \rho/\pa q=1$) to $0$ gives
\beq
\frac{\pa\rho}{\pa q}(t,x)= e^{G(s_2)}
\eq
Since $\rho$ is constant along the integral curves $X(s)$ it follows from
\eqref{eq:HLLassumption1}
\beq
|G(s_2)|\leq \frac{c_1\varepsilon}{2(1+|\rho|)^\nu}
\int_{s_2}^0\frac{ds}{1+X^0(s)}
\leq
\frac{c_1\varepsilon}{(1+|\rho|)^\nu}
\int_{\tau_2}^t\frac{dt^\prime}{1+t^\prime}
\leq\frac{c_1\varepsilon}{(1+|\rho|)^\nu}\ln{\Big|\frac{1+t}{1\!+\!|\rho|}\Big|},
\eq
since $t=X^0$, $d X^0/ds= L_2^0=1-\frac{1}{4}H_{LL}\geq 1/2$
and $\tau_2=t_2=2|\rho|$.
This proves \eqref{eq:eikonal1}.

 \eqref{eq:eikonal2} follows as above from integrating
\eqref{eq:qchar} in the form
\beq
\Big|\frac{d}{d s} \ln{\big|1\!+\!|\, q|\big|}\Big|=\frac{1}{2}\frac{|H_{LL}|}{1\!+\!|\, q|}
\leq \frac{c_1\varepsilon}{1+t}.
\eq
\end{proof}

We can write \eqref{eq:rhoqeq} as
 \beq\label{eq:rhoqeq3}
L_2^\alpha\pa_\alpha \ln{|\,\rho_q|}=-\frac{1}{2}\pa_q H_{LL}.
\eq
Integrating along the integral curves from the boundary
where $|\,t-r|=t/2$ and $\rho_q=1$ gives
\beq\label{eq:rhoqeq4}
\ln{|\,\rho_q|}=-\frac{1}{2}\int_{2|\rho|-\rho/2}^{p/2}\pa_q H_{LL}\, ds
\eq

 We now give some further estimates for the approximate solution of the eikonal equation:

 \begin{lemma}\label{eikonaltwo} Suppose that the assumption of Proposition \ref{eikonalone}
 hold and
 \beq\label{eq:HLLassumption2}
 |\pa^2 H_{LL}|\leq \frac{c_2\, \varepsilon}{1+t} \Big|\frac{\pa \rho}{\pa q}\Big|
\frac{1}{(1+|\rho|)^{1+\nu}},
 \eq
 for some $\nu\geq 0$.
 Then with $\pa_\rho=\rho_q^{-1}\pa_q$ we have
 \beq\label{eq:eikonal6}
 |\pa_\rho\, \pa_q\rho|\le \frac{c_2\varepsilon}{
 (1+|\rho|)^{1+\nu}} \pa_q \rho\, \ln{\Big|\frac{1+t}{1\!+\!|\rho|}\Big|}.
 \eq
 \end{lemma}
 \begin{proof}
Since $[L_2^\alpha\pa_\alpha,\pa_\rho]=0$ we obtain by differentiating
\eqref{eq:rhoqeq3}:
\beq
 L_2^\alpha\pa_\alpha \pa_\rho \ln{|\,\rho_q|}
=-\frac{1}{2}\pa_\rho \pa_q H_{LL}
\eq
Hence by \eqref{eq:HLLassumption2},
\beq\label{eq:rhoqeq2}
\big|L_2^\alpha\pa_\alpha \pa_\rho \ln{|\,\rho_q|}\big|\leq
 \frac{c_2\,\varepsilon}{(1+t)\,(1+|\rho|)^{1+\nu}}.
\eq
Moreover by differentiating \eqref{eq:rhoqeq4} we see that
the initial condition on $|\,t-r|=t/2$ are
$\pa_q\ln{|\, \rho_q|}=(\operatorname{sign}{\rho}-1/4)\pa_q H_{LL}$
It therefore follows from integrating \eqref{eq:rhoqeq2} from
the boundary where $t=2|\rho|$:
\beq
\big|\pa_q \ln{|\, \rho_q|}\big|\leq
\frac{c_2\varepsilon}{
 (1+|\rho|)^{1+\nu}} \ln{\Big|\frac{1+t}{1\!+\!|\rho|}\Big|}
+\frac{c_1\varepsilon}{(1+t)(1+|\rho|)^\nu}
\eq
 Since $t\geq 2|\rho|$ in the domain where $\rho_q\neq 1$,
 \eqref{eq:eikonal6} follows from this.
 \end{proof}

 \begin{lemma}\label{eikonalthree} Suppose the assumption of Proposition \ref{eikonalone} hold
and
\beq\label{eq:eikonalassumptions2}
|H_{L\u L}|+|H_{LA}|+|H_{AB}|+|H_{LL}|+|\,\Omega H_{LL}|\leq \frac{c_2\varepsilon(1+|\rho|)^{1-\nu}}{(1+t)^{1-c_2^\prime\varepsilon}},
 \eq
 for some $\nu^{\,\prime\prime}\geq 0$.
 Then
 \beq\label{eq:tanderrho}
 |\pab \rho|
 \leq \frac{c_2^\prime \sigma}{1+\sigma} \,\, \frac{(1+|\rho|)^{1-\nu}}{(1+t)^{1-c_2^\prime\varepsilon}},
\qquad\quad \sigma
=c_1\varepsilon\ln{|1+t|},
 \eq
 and
 \beq\label{eq:eikonalrho}
 |g^{\alpha\beta}\rho_\alpha \rho_\beta-\delta^{AB}\pa_A \rho\,\, \pa_B\rho|\leq
c_2^\prime \varepsilon \frac{ (1+|\rho|)^{2-2\nu}}{(1+t)^{2-c_2^\prime\varepsilon}}.
 \eq
 \end{lemma}
 \begin{proof} We have $|\pab \rho|\les |\pa_p\, \rho|+\sum_{\,0<i<j}|\Omega_{ij} \rho|/(1+t)$,
where
\beq
\pa_p \,\rho=-\frac{H_{LL}}{4}\pa_q\rho,
\eq
and
\beq
L_2^\alpha\pa_\alpha \Omega\rho=-\frac{1}{2}(\Omega H_{LL})\, \, \pa_q \rho.
\eq
Using \eqref{eq:eikonalassumptions2} and integrating the last equation
from $t=2|\rho|$, where $\Omega_{ij}\rho=0$,
gives \eqref{eq:tanderrho}. (Here we wrote
$\int_0^t c\varepsilon (1+t)^{c\varepsilon-1}\, dt=(1+t)^{c\varepsilon}-1\leq e^{c\varepsilon \ln{|1+t|}}-1$
and used the inequality $e^a-1\leq Ce^a \, a/(1+a)$.)

 By \eqref{eq:eikonalframe}
 \beq
\big|g^{\alpha\beta}\pa_\alpha\rho\, \pa_\beta \rho-\delta^{AB}\pa_A \rho\,\, \pa_B\rho\big|\leq 2|\pa_q
\rho|\,| L_1^{\alpha}\pa_\alpha \rho|+C(|H_{LL}|+|H_{LA}|+|H_{AB}|)|\overline{\pa}\rho|^2.
 \eq
Here
\beq
|L_1^\alpha\pa_\alpha \rho|\leq |L_2^\alpha\pa_\alpha\rho|+\Big|\Big(\frac{1}{2}H_{L\u L} L^\alpha +H_{L}^{\,\, A} A^\alpha\Big) \pa_\alpha \rho\Big|\leq
\big(|H_{L\u L}|+|H_{LA}|\big) |\pab \rho|
\eq
Hence
 \beq
\big|g^{\alpha\beta}\pa_\alpha\rho\, \pa_\beta \rho
-\delta^{AB}\pa_A \rho\,\, \pa_B\rho\big|\les
\big(|H_{L\u L}|+|H_{LA}|\big)|\pa_q\rho|\, |\pab\rho| +(|H_{LL}|+|H_{LA}|+|H_{AB}|)|\overline{\pa}\rho|^2.
 \eq
 These estimates together with \eqref{eq:tanderrho} gives also
\eqref{eq:eikonalrho}.
 \end{proof}

\section{The sharp decay estimates for the nonlinear problem}
 In this section we start by assuming the weaker decay estimates
 \beq\label{eq:weakdecay}
 |Z^I \phi|\leq c_0\varepsilon(1+t)^{-\nu},\qquad |I|\leq N-3,\qquad c_0\varepsilon\leq 1,
 \eq
 for some $ 0<\nu<1$ and some sufficiently large $N$.
 We also assume that $\phi$ is a solution of the nonlinear problem
 with compactly supported data in the set $|x|\leq 1$, which means  that
 \beq\label{eq:support}
 \phi(t,x)=0,\qquad\text{when}\quad |x|\geq t+1, \quad\text{and}\quad t\geq 0.
 \eq
 \eqref{eq:weakdecay}
 can be obtained from energy estimates using the
 Klainerman-Sobolev inequalities. From the weak decay estimates we
 will derive stronger decay estimates. The stronger estimates will
 be derived in several steps. Since our metric
 $g^{\alpha\beta}=m^{\alpha\beta}+H^{\alpha\beta}$, where
 $H^{\alpha\beta}=H^{\alpha\beta}(\phi)$ are smooth functions of
 $\phi$ vanishing at the origin and by scaling we may assume that
 (so that \eqref{eq:weakdecay} holds also for $H^{\alpha\beta}$)
 \beq
\sum_{\alpha\beta} |H^{\alpha\beta}(\phi)|\leq \frac{1}{64}|\phi|
\eq
In what follows $C$ will denote universal constants that depend only on the
the particular functions $H^{\alpha\beta}(\phi)$, but are independent of $\phi$.
$c_0^\prime$ will denote a constant that is multiple of $c_0$ i.e. $C c_0$.
Constants of the form $c_k$ and $c_k^\prime=Cc_k$ depend only on $c_{k-1}$
and universal constants.
The estimates
\eqref{eq:weightedstrongdecay0}-\eqref{eq:weightedstrongdecay2} below, were used already in the spherically symmetric case in \cite{L1}.
\begin{prop}\label{strongdecay} Suppose that $\phi$ is a solution of the nonlinear equation for which  \eqref{eq:weakdecay} and \eqref{eq:support} hold.
Let $\rho=\rho_2$ be as in the previous section. Then there are constants
$c_1=Cc_0$ and $c_2=Cc_0$, for some universal constant $C$, independent of $\phi$ if $c_0\varepsilon\leq 1/C$, such that
\beq\label{eq:weightedstrongdecay0}
|\phi|\leq \begin{cases} c_1\varepsilon (1+t)^{-1}(1+|q|)\\
c_1\varepsilon (1+t)^{-1+c_1\varepsilon} (1+|\rho|)^{1-\nu-c_1\varepsilon}
\end{cases}
 \eq
\beq\label{eq:weightedstrongdecay1}
 |\pa\phi|\leq c_1\varepsilon (1+t)^{-1}(1+|\rho|)^{-\nu},
 \eq
 \beq\label{eq:weightedstrongdecay2}
 |\pa^2 \phi|\leq
 c_2\varepsilon(1+t)^{-1}(1+|\rho|)^{-1-\nu}|\pa \rho/\pa q|,
 \eq
and
\beq\label{eq:weightedstrongdecayZ}
|Z \phi|\leq \begin{cases} c_2\varepsilon (1+t)^{-1+c_2\varepsilon}(1+|\rho|)^{1-\nu}\\
 c_2 \varepsilon(1+t)^{-1} \big( |q|+(1+t)^{c_2\varepsilon}\big).
\end{cases}
 \eq
Moreover, there are constants $c_k$ depending only on $c_{k-1}$ such that
 \beq\label{eq:weightedstrongdecaykI}
|\pa^{\,\bold{k}} Z^I \phi|\leq
c_k\varepsilon(1+t)^{-1+c_k\varepsilon}(1+|q|)^{1-k-\nu},
 \qquad \max{(1,|\bold{k}|)}+|I|\leq N-4,
 \eq
where $k=|\bold{k}|$.
\end{prop}

 \subsection{The decay of the first order derivatives
 \eqref{eq:weightedstrongdecay0} and \eqref{eq:weightedstrongdecay1}}
 \label{decayfirstder}
 Since by \eqref{eq:weakdecay}
 condition (1) in Lemma \ref{strongdecaywaveeq} hold and it follows from
 \eqref{eq:weakdecay} that the right hand side of
 \eqref{eq:strongdecaywaveeq} is bounded so
 \beq \label{eq:strongdecay1}
 |\pa \phi|\leq c_0^\prime\varepsilon(1+t)^{-1}
 \eq
 The first estimate in \eqref{eq:weightedstrongdecay0} follows from integrating this from
 $r=t+1$ where $\phi=0$.
Hence
\beq\label{eq:firstHest}
|\pa H|+(1+|q|)^{-1} |H|\leq \frac{c_1\varepsilon}{1+t},
\eq
it follows that in
fact condition (2) in Lemma \ref{strongdecaywaveeq} also hold.
\eqref{eq:weightedstrongdecay1} therefore follows from Lemma \ref{weightedstrongdecaywaveeq} with $\mu=a=0$.
 The second estimate in \eqref{eq:weightedstrongdecay0} follows from integrating
\eqref{eq:weightedstrongdecay1} and using Proposition \ref{eikonalone},
which hold since we just showed that \eqref{eq:weightedstrongdecay1} hold.

\subsection{The sharp decay estimates for second order derivatives \eqref{eq:weightedstrongdecay2}}\label{decaysecondder}
 The estimate \eqref{eq:weightedstrongdecay2} essentially comes from that $\pa_\rho$ commutes with $L_2$. Since \eqref{eq:firstHest} hold:
\begin{lemma} Suppose that $H$ satisfy \eqref{eq:firstHest} and let
$\rho$ be as in the previous section. Then
\beq\label{eq:seconderestrho}
\Big|2L_2^\alpha\pa_\alpha( r \pa_\rho \pa_q \phi)-r\pa_\rho \Boxr_g\phi\Big|
 \leq\frac{C\rho_q^{-1}}{1+t+r}\sum_{|I|\leq 2} \,|\pa Z^I\phi|,
\qquad\text{where}\quad \pa_\rho =\rho_q^{-1}\pa_q.
\eq
and
\beq\label{eq:vectorfieldestrho}
\Big|2L_2^\alpha\pa_\alpha( r \pa_\rho \psi)-r\rho_q^{-1} \Boxr_g\psi\Big|
 \leq c_1\varepsilon |\pa_\rho\psi|+ \frac{C\rho_q^{-1}}{1+t+r}\sum_{|I|\leq 2} \,|Z^I \psi|,
\eq
\end{lemma}
\begin{proof}  Since $2L_2^\alpha\pa_\alpha \rho=(\pa_p+H_{LL} \pa_q )\rho=0$ it follows that
 $2L_2^\alpha\pa_\alpha \pa_q\rho=-\pa_q H_{LL}\, \, \pa_q\rho$ and
 $2L_2^\alpha\pa_\alpha \rho_q^{-1}=\rho_q^{-1} \pa_q H_{LL}$.
Hence
\beq
2L_2^\alpha\pa_\alpha( r \rho_q^{-1}\pa_q^2 \phi)=\rho_q^{-1} \big(2L_2^\alpha\pa_\alpha (r\pa_q^2\phi)+(r\pa_q H_{LL})\pa_q^2\phi\big)
\eq
and \eqref{eq:seconderestrho} follows from \eqref{eq:secondereq}.
\eqref{eq:vectorfieldestrho} follows in the same way.
\end{proof}
It follows from \eqref{eq:seconderestrho} using \eqref{eq:eikonal1} and \eqref{eq:eikonal2} that
 \beq
 \Big|L_2^\alpha\pa_\alpha\Big(\rho_q^{-1}
r\,\pa_q^2\phi\Big)\Big|\leq
\frac{C\rho_q^{-1}}{1\!+t} \sum_{|I|\leq 2\!\!}|\pa Z^I\phi|
\leq \frac{C\rho_q^{-1}}{(1\!+t)(1+|q|)} \sum_{|I|\leq 3\!\!}|Z^I\phi|
\leq \frac{c_0^\prime\varepsilon}{(1+t)^{1+\nu-2c_1\varepsilon} (1+|\rho|)^{1+2c_1\varepsilon}}.
 \eq
If we as in the proof of Lemma \ref{strongdecaywaveeq} integrate from $r=t/2$, where $t\sim\rho$
and $|\rho_q^{-1} r\pa_q^2\phi|\leq c_0 \varepsilon (1+|\rho|)^{-\nu}$ we get
\beq
|\rho_q^{-1}
r\pa_q^2\phi|\leq c_1^\prime \varepsilon (1+|\rho|)^{-1-\nu},
\eq
since we assumed that $c_0 \varepsilon\leq 1$.
\eqref{eq:weightedstrongdecay2} follows from this
using Lemma \ref{tanderZ} and \eqref{eq:weakdecay}.

\subsection{The decay estimate for one vector field \eqref{eq:weightedstrongdecayZ}}\label{decayonevectorfield}
Since $\Boxr_g\phi=0$ we have by \eqref{eq:onecommute}:
 \beq
 |\Boxr_g Z\phi|=
 \big|\big( Z H^{\alpha\beta}+2 C_{Z\gamma}^{\,\alpha}  H^{\gamma\beta}+C_Z H^{\alpha\beta}\big)\pa_\alpha\pa_\beta\phi\big|\leq C\big(|Z\phi|+|\phi|\big)\,|\pa^2\phi|,
 \eq
 and hence by \eqref{eq:vectorfieldestrho} applied to $\psi=Z\phi$;
 \beq
 |L_2^\alpha\pa_\alpha \big(r\pa_\rho Z \phi\big)|\leq
 C r(|Z\phi|+|\phi|)\frac{|\pa^2\phi|}{\rho_q}
 +c_1\varepsilon|\pa_\rho Z\phi|+\frac{C\rho_q^{-1}}{1+t} \sum_{|J|\leq 3} |Z^J\phi|.
 \eq
Hence using \eqref{eq:weightedstrongdecay2}, \eqref{eq:weightedstrongdecay1},
\eqref{eq:weakdecay} and \eqref{eq:weightedstrongdecay0} we get
\beq
 \Big|L_2^\alpha\pa_\alpha \big(r\, \pa_\rho Z\phi\big)\Big|\leq
c_2^\prime\varepsilon\Big( \frac{|Z\phi|}{1+|\rho|} + |\pa_\rho Z\phi|\Big) +
\frac{c_0^\prime\varepsilon(1\!+\! |\, \rho|)^{-c_1\varepsilon} }{(1+t)^{1+\nu-c_1\varepsilon}}
+\frac{c_2^\prime \varepsilon^2 }{(1+t)^{1-c_1\varepsilon}(1+|\rho|)^{2\nu+c_1\varepsilon}}.
 \eq
Since also
\beq
|Z \phi|  +(1+|\rho|)|\pa_\rho Z\phi|=
|Z\phi|+(1+|q|) |\pa Z\phi|
\les \sum_{|I|\leq 2} |Z^I\phi|\les \varepsilon
(1+t)^{-\nu},\qquad\text{when $|\,t\!-r|=t/2$},
\eq
it follows from Lemma \ref{Zlemma} below that
\beq
|Z \phi| (1+|\rho|)^{-1} +|\pa_\rho Z\phi|\leq c\varepsilon
(1+t)^{-1+c\varepsilon} (1+|\rho|)^{-\nu},\qquad
c=16(c_2^\prime+c_0^\prime)
\eq
The desired inequality \eqref{eq:weightedstrongdecayZ} follows from this and Lemma \ref{convexity}, since $(1+|\rho|)^{c_1\varepsilon}\les
(1+t)^{c_1\varepsilon}$.

\begin{lemma}\label{Zlemma} Suppose that for some $\nu>0$ we have
\beq\label{eq:psirhoder}
(1+|\rho|)|\pa_\rho\psi|+|\psi|\leq c_0^\prime(1+|\rho|)^{-\nu}\!\!\!,
\,\,\, \text{when $|\,t\!-r|=t/2$ or $t+r\leq 2$.}
\eq
and
\beq\label{eq:supportpsi}
\psi=0,\,\,\text{when $r\!>\!t+1$ and $t\!>\!0$.}
\eq
 Suppose also that
 \beq
\Big|L_2^\alpha\pa_\alpha \big(r\, \pa_\rho\psi \big)\Big|\leq
c_2^\prime \varepsilon\Big(\frac{|\psi|}{1+|\rho|}+|\pa_\rho \psi|\Big) +
\frac{c_0^\prime\varepsilon(1\!+\! |\, \rho|)^{-c_1\varepsilon} }{(1+t)^{1+\nu-c_1\varepsilon}}
+\frac{\gamma\varepsilon^2 }{(1+t)^{1-\gamma\varepsilon}(1+|\rho|)^{\nu+\gamma\varepsilon}}.
 \eq
Then
\beq
|\psi| (1+|\rho|)^{-1} +|\pa_\rho \psi|\leq c\varepsilon
(1+t)^{-1+c\varepsilon} (1+|\rho|)^{-\nu},\qquad
c=16(c_2^\prime+c_0^\prime+\gamma)
\eq
 \end{lemma}
\begin{proof}
 If we now introduce the new variables $p=t+r=2s$, $q=r-t=q(\rho,s,\omega)$
 given by \eqref{eq:changeofvariables}:
 $\widetilde{v}(\rho,s,\omega)=v(q,p,\omega)$.
 Then
 \beq
L_2^\alpha\pa_\alpha v(q,p,\omega)=\big[\pa_s \widetilde{v}(\rho,s,\omega)\big]\Big|_{\rho=\rho(q,p,\omega),\, p=2s}.
\eq

It is also easy to see that if we substitute $r=(p+q)/2$ in the left of
then the term with $q/2$ in place of $r$ in the left can be bounded by terms of the form already included in the right so with $\psi(\rho,s,\omega)=(Z^I\phi)\big(q(\rho,s,\omega),2s,\omega\big)$ we have
\beq
 \Big|\pa_s \big( s\, \pa_\rho \psi\big)\Big|\leq
 c_2^\prime\varepsilon \Big(\frac{|\psi|}{1+|\rho|} + |\pa_\rho \psi|\Big)+
\frac{c_0^\prime\varepsilon(1\!+\! |\, \rho|)^{-c_1\varepsilon} }{(1+s)^{1+\nu-c_1\varepsilon}}
+\frac{\gamma\varepsilon^2 }{(1+s)^{1-\gamma\varepsilon}(1+|\rho|)^{\nu+\gamma\varepsilon}}
 \eq
If we integrate this from the boundary of
$D=\{(r,t);\, t/2< r<3t/2\}=\{(\rho,p);\, -2s/3<\rho<2s/5\}$ (since $q\!=\!\rho$
and $s\!=\!3\rho/2$ or $s\!=\!5\rho/2$ on the boundary)  using the bound \eqref{eq:psirhoder} on the boundary we get
\beq
 s\, (1+|\rho|)^{\nu}\, |\pa_\rho\psi| \leq  c_2^\prime \varepsilon
 \int_{c_\pm |\rho|}^{s}\Big(
\frac{|\psi|}{1+|\rho|} + |\pa_\rho \psi|\Big)\, d s
+c_0^\prime \varepsilon+\varepsilon(1+p)^{\gamma\varepsilon}
(1\!+\! |\,\rho|)^{-\gamma\varepsilon}
\eq

For any $0\leq a<1$ we have
\beq
|\psi(\rho,s,\omega)-\psi(1,s,\omega)|\leq \int^{1}_\rho \!\! |\pa_\rho \psi(\rho,s,\omega)|\, d\rho
\leq (1+|\rho|)^{1-a}\!\!\! \sup_{\rho\,\leq \varrho\,\,\leq 1}\! (1+|\varrho|)^{a}|\pa_\rho \psi(\varrho,s,\omega)|.
\eq
By \eqref{eq:supportpsi} $\psi(1,s,\omega)=0$, when $s\geq 1$ and
by \eqref{eq:psirhoder} $|\psi(1,s,\omega)|+|\pa_\rho \psi(1,s,\omega)|\leq c_0\varepsilon$ for all $|s|\leq 1$.
With
\beq
M(s)=\sup_{\rho\,\in D_s} (1+|s|)\, (1+|\rho|)^{\nu}\, |\pa_\rho\psi(\rho,s,\omega)|,
\qquad\text{where}\quad  D_s=\{\rho; (s,\rho,\omega)\in D\}.
\eq
we hence have
\beq
M(s)\leq 2 c_2^\prime\,\varepsilon \int_1^s \frac{M(s^{\,\prime})}{1+|s^{\,\prime}|}\, ds^{\,\prime}+c_0^\prime\varepsilon+\varepsilon(1+|s|)^{\gamma\varepsilon}
\eq
If $G(s)$ denotes the integral we hence have $(1+|s|)G^{\,\prime}(s)= M(s)\leq
2c_2^\prime \varepsilon G(s)+c_0^\prime\varepsilon+\varepsilon (1+|s|)^{\gamma\varepsilon}(1\!+\! |\,\rho|)^{-\gamma\varepsilon}$,
and if we multiply by the integrating factor $(1+|s|)^{-c\varepsilon} $ we get
\beq
d \big(G(s) (1+|s|)^{-c\varepsilon}\big)/d s \leq  (1+|s|)^{-c\varepsilon-1}
\big(c_0^\prime\varepsilon+\varepsilon(1+|s|)^{\gamma\varepsilon}\big).
\eq
If we integrate this from $1$ to $s$ we get $G(s)\leq 4(1+|s|)^{c\varepsilon}$,
and hence
\beq
M(s)\leq 2^{-1} c\varepsilon (1+|s|)^{c\varepsilon}.
\eq
We conclude that
\beq
|\psi| (1+|\rho|)^{-1} +|\pa_\rho \psi|\leq c\varepsilon
(1+t)^{-1+c\varepsilon} (1+|\rho|)^{-\nu}.
\eq
\end{proof}

\begin{lemma}\label{convexity} If $\nu>0$ then
\beq
 (1+t)^{c\varepsilon}
(1+|\, q|)^{1-\nu}
\leq (1+t)^{c\varepsilon/\nu} +(1+|q|)
\eq
\end{lemma}
\begin{proof}
The follows from the inequality
\beq
a^{\nu} b^{1-\nu} \leq a+b ,\qquad 0<\nu<1, \qquad a>0,\quad b>0.
\eq
\end{proof}

 \subsection{The decay estimates for higher order derivatives}
\label{sharpdecayfour}
 Let us first prove that:
 \begin{lemma} For $1\leq k\leq N-4$ we have
\beq \label{eq:weightedstrongdecayk}
 |\pa^{\,\bold{k}} \phi|\leq
\frac{c_k\varepsilon}{1+t}(1+|\rho|)^{1-k-\nu}
\Big(\frac{1+t}{1\!+\!|\,\rho|}\Big)^{c_k\varepsilon V(\rho)} ,\qquad\quad V(\rho)=(1\!+\! |\,\rho|)^{-\nu},
 \quad |\bold{k}|=k,
 \eq
 if $\varepsilon$ is sufficiently small.
\end{lemma}
\begin{proof} We will prove the lemma by induction.
 If $k=1$ we already proved a stronger estimate
 and differentiating the equation $\Boxr_g \phi=m^{\alpha\beta}\pa_\alpha\pa_\beta\phi+H^{\alpha\beta}\pa_\alpha\pa_\beta\phi=0$
gives
 \beq
\Boxr_g \pa^{\,\bold{n}} \phi= -\sum_{\bold{m}+\bold{k}=\bold{n}, \,
|\bold{m}|\geq 1} \pa^\bold{m}H^{\alpha\beta}\,\,
\pa_\alpha\pa_\beta \pa^\bold{k} \phi.
 \eq
 Hence by \eqref{eq:Hder}
 \beq\label{eq:boxmanyder}
 |\Boxr_g \pa^{\,\bold{n}} \phi|\leq  C |\pa \phi| \sum_{|\bold{m}|
 =|\bold{n}|+1} |\pa^\bold{m} \phi|
 +C\sum_{|\bold{k}_1|+\cdots+|\bold{k}_\ell|=|\bold{n}|+2,\,
 1\leq |\bold{k}_j|\leq |\bold{n}|, \, \ell\geq 2}
 |\pa^\bold{k_1} \phi|\cdots |\pa^{\bold{k}_{\,\ell}} \phi|.
 \eq
 Using \eqref{eq:weightedstrongdecay1} and \eqref{eq:weightedstrongdecayk} for $|\bold{k}|\leq n=|\bold{n}|$ (and the fact that $|\rho|\les t$), we hence obtain
\beq
 |\Boxr_g \pa^{\,\bold{n}} \phi|\leq \frac{c_1^\prime \varepsilon}{(1+t)(1+|\rho|)^\nu}
  \sum_{|\bold{m}|
 =n+1} |\pa^\bold{m} \phi|
+\frac{(c_n^\prime)^2\varepsilon^2}
 {(1+t)^2 (1+|\rho|)^{n+2\nu}}\Big(\frac{1+t}{1\!+\!|\,\rho|}\Big)^{c_n^\prime\varepsilon V(\rho)} ,
 \eq
if $\varepsilon$ is so small that $c_n\varepsilon\leq 1$.
 By \eqref{eq:derframeZ}, \eqref{eq:eikonal2} and \eqref{eq:weakdecay}
 \beq\label{eq:Zq}
 |Z^I \pa^\bold{n}\phi|\leq C (1+|q|)^{-n} \!\!\!\!\!\sum_{|J|\leq n+|I|}\!\!\! |Z^J\phi|
\leq c_0^\prime \varepsilon (1+|q|)^{-n} (1+t)^{-\nu}\leq
c_0^\prime \varepsilon (1+t)^{-\nu+n c_1\varepsilon}(1+|\rho|)^{-n-n c_1\varepsilon}.
 \eq
The lemma now follows from the following lemma:
\end{proof}
\begin{lemma} Suppose that
\beq\label{eq:psiest}
\sum_{|\bold{n}|=n,\, |I|\leq 2}|Z^I\psi_\bold{n}|\leq
c_0^\prime \varepsilon (1+t)^{-\nu+n c_1\varepsilon}(1+|\rho|)^{-n-n c_1\varepsilon}
\eq
and
\beq
\sum_{|\bold{n}|=n} |\Boxr_g \psi_\bold{n}|\leq \frac{c_1^\prime \varepsilon}{(1+t)(1+|\rho|)^\nu}
  \sum_{|\bold{n}|
 =n} |\pa\psi_\bold{n}|
+\frac{(c_n^\prime)^2\varepsilon^2}
 {(1+t)^2 (1+|\rho|)^{n+2\nu}}\Big(\frac{1+t}{1\!+\!|\,\rho|}\Big)^{c_n^\prime\varepsilon V(\rho)}.
\eq
Then
\beq
\sum_{|\bold{n}|=n}|\pa \psi_\bold{n}|\leq
\frac{c_{n+1}\varepsilon}{1+t}
\Big(\frac{1+t}{1\!+\!|\,\rho|}\Big)^{c_{n+1}\varepsilon V(\rho)}(1+|\rho|)^{-n-\nu}
\eq
\end{lemma}
\begin{proof}
Using \eqref{eq:phiwave2} as in the proof of Lemma \ref{strongdecaywaveeq} we have
\beq
2
|L_2^{\alpha}\pa_\alpha(r\pa_q \psi_\bold{n})|\leq r|\Boxr_g \psi_\bold{n}| +\frac{C}{1+t+r}\!\sum_{|I|\leq 2}
|Z^I\psi_\bold{n}|,
\eq
 Using that
 \beq\label{eq:derqZ}
  |\pa \psi_\bold{n}|\leq C|\pa_q \psi_\bold{n}| + C (1+t+r)^{-1} \sum_{|I|=1} |Z^I  \psi_\bold{n}|
\eq
we have
\beq
|L_2^{\alpha}\pa_\alpha(r\pa_q \psi_\bold{n})|\les
\frac{c_1^\prime \varepsilon \sum_{|\bold{m}|=|\bold{n}|}|r\pa_q\psi_\bold{m}|}{(1+t)(1+|\rho|)^\nu}
+\frac{(c_n^\prime)^2\varepsilon^2}
 {(1+t) (1+|\rho|)^{n+2\nu}}\Big(\frac{1+t}{1\!+\!|\,\rho|}\Big)^{c_n^\prime\varepsilon V(\rho)}
 +\frac{c_0^\prime\varepsilon(1+|\rho|)^{-n-n c_1\varepsilon}}{(1+t)^{1+\nu-n c_1\varepsilon}}.
\eq
 Since
 $\big|L_2^\alpha\pa_\alpha \sum_{|\bold{n}|=n}|r\pa_q\psi_\bold{n}|\big|\leq
\sum_{|\bold{n}|=n}|L_2^\alpha\pa_\alpha(r\pa_q\psi_\bold{n})|$, and $L_2^\alpha\pa_\alpha\rho=0$,
we have
\beq
L_2^\alpha\pa_\alpha M_n\leq \frac{c_1^\prime\varepsilon M_n}{(1+t)(1+|\rho|)^\nu}
+\frac{(c_n^\prime)^2 \varepsilon^2}
 {(1+t) (1+|\rho|)^{\nu}}\Big(\frac{1+t}{1\!+\!|\,\rho|}\Big)^{c_n^\prime\varepsilon V(\rho)}
 +\frac{c_0^\prime\varepsilon(1+|\rho|)^{\nu-n c_1\varepsilon}}{(1+t)^{1+\nu-n c_1\varepsilon}},
\eq
where
\beq
M_n=\sum_{|\bold{n}|=n}|r\pa_q\psi_\bold{n}|(1+|\rho|)^{n+\nu}.
\eq
Let $c_{n+1}=2(c_1+c_{n}^\prime+1+4n c_0)$ and
\beq
N_n=\sum_{|\bold{n}|=n}|r\pa_q \psi_\bold{n}|(1+|\rho|)^{n+\nu}
\Big(\frac{1+t}{1\!+\!|\,\rho|}\Big)^{-c_{n+1}\varepsilon V(\rho)}.
\eq
Then
\beq
L_2^\alpha\pa_\alpha N_n\leq \frac{(c_n^\prime)^2\varepsilon^2}
 {(1+t) (1+|\rho|)^{\nu}}
\Big(\frac{1+t}{1\!+\!|\,\rho|}\Big)^{-\varepsilon V(\rho)}
 +\frac{c_0^\prime\varepsilon(1+|\rho|)^{\nu-n c_1\varepsilon}}{(1+t)^{1+\nu-n c_1\varepsilon}},
\eq
If we integrate along the integral curves of the vector field $L_2^\alpha\pa_\alpha$ from a point $(t_2,x_2)$ with $|t_2-|x_2||=t_2/2=\rho$
to a point $(t,x)$ as in the proof of Lemma \ref{strongdecaywaveeq}
we get
\beq
N_n(t,x)\leq N_n(t_2,x_2) + \int_{2\rho}^t\Big(
\frac{(c_n^\prime)^2 \varepsilon^2}
 {(1+t) (1+|\rho|)^{\nu}}
\Big(\frac{1+t}{1\!+\!|\,\rho|}\Big)^{-\varepsilon V(\rho)}
 +\frac{c_0^\prime \varepsilon(1+|\rho|)^{\nu-n c_1\varepsilon}}{(1+t)^{1+\nu-n c_1\varepsilon}} \Big) dt
\leq (c_n^\prime)^2\varepsilon,
\eq
since by \eqref{eq:psiest} $N_n(t,x)\leq c_0^\prime \varepsilon$ when $|t-|x||=t/2$.
The lemma now follows from the bound for $N_n$, \eqref{eq:derqZ}
and \eqref{eq:psiest}.
\end{proof}

\subsection{The decay estimates for more vector fields}
\label{strongdecaymanyvectorfields}
We will use induction to prove that
\beq\label{eq:weightedstrongdecaykIcopy}
|\pa^{\,\bold{k}} Z^I \phi|\leq
c_{k,i}\,\varepsilon(1+t)^{-1+c_{k,i}\varepsilon}(1+|\rho|)^{1-k-\nu}\!\!\!,
 \qquad \max{(1,k)}+i\leq N\!-4,\quad k\!=\!|\bold{k}|,\,\,\,\, i\!=\!|I|.
 \eq
Note that by \eqref{eq:eikonal2} $(1+|q|)( 1+t)^{-c_1\varepsilon}\leq (1+|\rho|)\leq (1+|q|)( 1+t)^{c_1\varepsilon}$
so we could just as well have stated \eqref{eq:weightedstrongdecaykIcopy} with
$\rho$ replaced by $q$.

We will use induction in $|I|$,
and for fixed $|I|$ induction in $|\bold{k}|$.
We will start by proving \eqref{eq:weightedstrongdecaykIcopy} for $|I|=0$ and
all $|\bold{k}|$.
Then we prove \eqref{eq:weightedstrongdecaykIcopy} for $|I|=m\geq 1$ and $|\bold{k}|\leq 1$ assuming  \eqref{eq:weightedstrongdecaykIcopy} for
$|I|\leq m-1$ and all $|\bold{k}|$.
Finally we prove \eqref{eq:weightedstrongdecaykIcopy} for $|I|=m$ and $|\bold{k}|=n+1\geq 2$ assuming \eqref{eq:weightedstrongdecaykIcopy} for $|I|=m$ and $|\bold{k}|\leq n$ and \eqref{eq:weightedstrongdecaykIcopy} for $|I|\leq m-1$ and all $|\bold{k}|$.

{\bf Proof of \eqref{eq:weightedstrongdecaykIcopy} for $|I|=0$ and all $|\bold{k}|$.}
In \eqref{eq:weightedstrongdecayk} we have already proven a stronger estimate than \eqref{eq:weightedstrongdecaykIcopy}
for $|I|=0$ apart from the case of $|\bold{k}|=0$
which follows from
integrating the same estimate for $|\bold{k}|=1$ in the $t-r$ direction,
using that $\phi$ vanishes when $r-t\geq 1$ and $t>0$.

{\bf Proof of \eqref{eq:weightedstrongdecaykIcopy} for
$|I|=m\geq 1$ and $|\bold{k}|\leq 1$
assuming \eqref{eq:weightedstrongdecaykIcopy} for
$|I|\leq m-1$ and all $|\bold{k}|$.}
By \eqref{eq:manycommute} and \eqref{eq:Hder} and the fact that $|Z^J\phi|\leq 1$ by \eqref{eq:weakdecay} we have
 \beq
 |\Boxr_g Z^I\phi|\leq C |Z^I \phi|\, |\pa^2\phi|\,\,
+C\sum_{|J|+|K|\leq |I|,\,\,  |J|\leq |I|-1,\,\,  |K|\leq |I|-1
\!\!\!\!\!\!\!\!\!\!\!\!\!\!\!\!\!\!\!\!\!\!\!\!\!\!\!\!\!\!
\!\!\!\!\!\!\!\!\!\!\!\!\!}
 |Z^{J}\phi|\, |\pa^2 Z^{K}  \phi|
 \eq
 and hence by \eqref{eq:vectorfieldestrho} applied to $\psi=Z^I\phi$;
 \beq
 |L_2^\alpha\pa_\alpha \big(r\pa_\rho Z^I\phi\big)|\leq
 C r|Z^I \phi|\frac{|\pa^2\phi|}{\rho_q} +c_1^\prime\varepsilon|\pa_\rho Z^I\phi|
 +\frac{C\rho_q^{-1}}{1+t} \sum_{|J|\leq \,|I|+2\!\!\!\!\!\!\!\!\!\!\!\!\!\!\!\!\!\!\!} |Z^J\phi|
+ \frac{C}{\rho_q}\sum_{|J|, |K| \leq |I|-1\!\!\!\!\!\!\!\!\!\!\!\!\!
\!\!\!\!\!\!\!\!\!}
 r\, |Z^J \phi|\, |\pa^2 Z^K \phi|
 \eq
Hence using \eqref{eq:weightedstrongdecay2}, \eqref{eq:weightedstrongdecay1}
and \eqref{eq:weightedstrongdecaykIcopy} for $|I|$ replaced by $|I|-1$ we get
\beq
 \Big|L_2^\alpha\pa_\alpha \big(r\, \pa_\rho Z^I\phi\big)\Big|\leq
 c_2^\prime\varepsilon\Big(\frac{|Z^I\phi|}{1+|\rho|} +|\pa_\rho Z^I\phi|\Big)+
\frac{c_0^\prime\varepsilon (1\!+\! |\, \rho|)^{-c_1\varepsilon}}{(1+t)^{1+\nu-c_1\varepsilon}}
+\frac{c_{0,m-1}^\prime c_{2,m-1}^\prime\varepsilon^2 }{(1+t)^{1-c\varepsilon}(1+|\rho|)^{2\nu}}
 \eq
It follows from Lemma \ref{Zlemma} that with
$c=16(c_2^\prime+c_0^\prime+c_{0,m-1}^\prime c_{2,m-1}^\prime)$ we have
\beq
|Z^I\phi| (1+|\rho|)^{-1} +|\pa_\rho Z^I \psi|\leq c\varepsilon
(1+t)^{-1+c\varepsilon} (1+|\rho|)^{-\nu}
\eq

{\bf Proof of \eqref{eq:weightedstrongdecaykIcopy} for
$|I|=m\geq 1$ and $|\bold{k}|=n+1\geq 2$
assuming \eqref{eq:weightedstrongdecaykIcopy} for
$|I|\leq m$ and all $|\bold{k}|\leq n$
and \eqref{eq:weightedstrongdecaykIcopy} for
$|I|\leq m-1$ and all $|\bold{k}|$.}

It follows from \eqref{eq:manycommutemixed} and \eqref{eq:Hder} that
 \beq
 |\Boxr_g \pa^{\bold{n}} Z^I\phi|\les \sum_{|\bold{m}|=|\bold{n}|+1
 \!\!\!\!\!\!\!\!\!\!\!\!\!\!\!\!} |\pa \phi|\, |\pa^{\bold{m}} Z^I\phi|
+\sum_{|{\bold{k}_1}|+\cdots+|\bold{k}_\ell|+|\bold{m}|=|\bold{n}|+2,\,
\, |\bold{k}_j|\leq |\bold{n}|,\, \ell\geq 1,\,\,
 |J_1|+\cdots+|J_\ell|+|K|\leq |I|,\,\,  |\bold{m}|\leq |\bold{n}|\text{ or }|K|<|I|
\!\!\!\!\!\!\!\!\!\!\!\!\!\!\!\!\!\!\!\!\!\!\!\!\!\!\!\!\!\!\!\!\!\!\!\!\!\!\!\!
\!\!\!\!\!\!\!\!\!\!\!\!\!\!\!\!\!\!\!\!\!\!\!\!\!\!\!\!\!\!\!\!\!\!\!\!\!\!\!\!
\!\!\!\!\!\!\!\!\!\!\!\!\!\!\!\!\!\!\!\!\!\!\!\!\!\!\!\!\!\!\!\!\!\!\!\!\!\!\!\!} |\pa^{\bold{k}_1} Z^{J_1}\phi|\cdots |\pa^{\bold{k}_\ell} Z^{J_\ell}\phi|\, \, |\pa^{\bold{m}}  Z^K \phi|.
 \eq
 Using \eqref{eq:weightedstrongdecay1} and \eqref{eq:weightedstrongdecaykIcopy} for $|\bold{k}|\leq n=|\bold{n}|$ and $|I|\leq m$, and
$|\bold{k}|\leq n+2$ and $|I|\leq m-1$ we hence obtain
 \beq
 |\Boxr_g \pa^{\bold{n}} Z^I\phi|\leq
 \frac{c_1\varepsilon}{(1+t)(1+|\rho|)^\nu}  \sum_{|\bold{m}|=|\bold{n}|+1} |\pa^\bold{m} Z^I\phi|
+\frac{(c_n^\prime)^2\varepsilon^2 }{(1+t)^{2-c\varepsilon}(1+|\rho|)^{|\bold{n}|+2\nu}}.
\eq
\eqref{eq:weightedstrongdecaykIcopy} for $|I|=m$ and $|\bold{k}|=n+1$ now follows
as in the proof of \eqref{eq:weightedstrongdecayk}.

\section{Weighted Energy estimates for the wave equation on a
curved background}

  We now
establish the basic energy identities with weight for solutions of
the equation \beq \label{eq:quasinh} \Boxr_{g} \phi =F \eq The
weight will be of the form
\beq\label{eq:energyweight}
w=e^{\,\sigma V(\rho)},\quad \sigma=\kappa\,\varepsilon\ln{|1+t|},   \quad V(\rho)=
|\rho-2|^{-\nup}, \quad \rho \leq 1,\qquad \nu^\prime,\kappa\geq 0
\eq
We note that by \eqref{eq:eikonal6}:
\beq
 g^{\alpha\beta}\rho_\alpha \rho_\beta\geq \delta^{AB}\pa_A \rho\,\, \pa_B\rho
-\big| g^{\alpha\beta}\rho_\alpha \rho_\beta-\delta^{AB}\pa_A \rho\,\, \pa_B\rho\big|
\geq -c_2^\prime \varepsilon \frac{ (1+|\rho|)^{2-2\nu^{\,\prime\prime}}}{(1+t)^{2-c_2^\prime\varepsilon}},
\eq
so $\rho$ satisfies the assumption below
if $c_2^\prime\varepsilon\leq 1/(\kappa\nup)$.
The following lemma was essentially proven in
\cite{A2}:
\begin{lemma}\label{curvedenergyest}
\label{lemma:Energy} Let $\phi$ be a solution of the equation
\eqref{eq:quasinh} decaying sufficiently fast as $|x|\to\infty$,
with a metric $g$ and a weight function $w$ as in \eqref{eq:energyweight}, with $\rho$, satisfying the conditions
 \beq
\label{eq:esmall} |g-m|\le \frac 12,\qquad |\pa g|\leq
\frac{c_1\varepsilon}{1+t}, \qquad
\rho_t<0,\qquad
 \frac{g^{\alpha\beta}\rho_\alpha
\rho_\beta}{\rho_t\,(1+|\rho|)^{1+\nup}} \geq -\frac{1/(\kappa\nup)}{(1+t)\ln{|1+t\,|}}
\eq
and $g^{\alpha\beta}=m^{\alpha\beta}$, the Minkowski metric, when $r>t+1$.
 Then for functions $\phi$ vanishing for $r>t+1$ we have, with $c=c_1+\kappa$;
 \beq
 \int_{\Si_{t}}\!\!|\pa\phi|^{2} w\, dx\leq
 4\! \int_{\Si_{0}} \!\!\! |\pa\phi|^{2} w\, dx
 +\int_0^t\!\!\frac{\,\,4 c\varepsilon\!}{1\!+\tau}\!
 \int_{\Sigma_\tau} \!\!\! |\pa\phi|^2\,w\,  dx \, d\tau
 +\frac{4\!}{c\varepsilon\,}\!\!\int_0^t\!\!
 \int_{\Sigma_\tau}\!\!\!
 (1+\tau)|\Boxr_g \phi|^2\,w\,  dx \, d\tau.
\eq
\end{lemma}
\begin{proof} Let $\phi_i=\pa_i \phi$, $i=1,2,3$, and $\phi_t=\pa_t\phi$.
If we differentiate below the integral sign and integrate by parts
we get
\begin{multline}
\frac{d}{dt} \int\big(-g^{00}\phi_t^2+g^{ij}\phi_i\phi_j\big) w\,
dx-\int 2\pa_j\big(g^{0j}\phi_t^2 w\big)\, dx
=\int 2\big(-g^{00} \phi_t\phi_{tt}+g^{ij}\phi_i\phi_{tj}-2g^{0j}\phi_t\phi_{tj}\big)
w\, dx \\
+\,\int  \big(-(\pa_t g^{00}) \phi_t^2+(\pa_t
g^{ij})\phi_i\phi_{j}-2(\pa_j g^{0j})\phi_t^2\big)w\, +
\big(-g^{00}\phi_t^2 w_t+g^{ij}\phi_i\phi_j w_t
-2g^{0j}\phi_t^2 w_j\big)\, dx\\
=\int 2\big(-g^{00}
\phi_t\phi_{tt}-g^{ij}\phi_t\phi_{ij}-2g^{0j}\phi_t\phi_{tj}\big)
 w \, dx \\
+\,\int  \big(-(\pa_t g^{00}) \phi_t^2+(\pa_t g^{ij})\phi_i\phi_{j}
-2(\pa_j g^{0j})\phi_t^2-2(\pa_i g^{ij})\phi_t\phi_j\big)w\, dx\\
+\int \big(-g^{00}\phi_t^2w_t+g^{ij}\phi_i\phi_j w_t-2
g^{0j}\phi_t^2 w_j -2g^{ij}\phi_t \phi_j w_i\big)\, dx
\end{multline}
Hence, since we also have assume that $\phi_t$ and $g^{0j}$ decay
fast enough that the boundary term vanishes at infinity
\begin{multline}
\frac{d}{dt} \int\big(-g^{00}\phi_t^2+g^{ij}\phi_i\phi_j\big) w\, dx
=\int w\big( \phi_t \,\Boxr_g \phi+(\pa_t g^{\alpha\beta})
\phi_\alpha \phi_\beta
-2(\pa_\alpha g^{\alpha\beta})\phi_\beta \phi_t\big)\, dx \\
+\int g^{\alpha\beta}\phi_\alpha\phi_\beta\, w_t -2\phi_t
g^{\alpha\beta}\phi_\alpha w_\beta\, dx
\end{multline}
Now
 \beq\label{eq:wder}
w_t=\frac{\kappa\nup\varepsilon\ln{|1+t|}}{|\rho-2|^{1+\nup}} \rho_t w
+\frac{\kappa\varepsilon}{(1+t)|\rho-2|^\nup} w,\qquad
w_i= \frac{\kappa\nup\varepsilon\ln{|1+t|}}{|\rho-2|^{1+\nup}} \rho_i w
 \eq
 If we set $\widehat{\phi}_\alpha=\phi_\alpha/\phi_t$ and
$\widehat{\, \rho}_\alpha=\rho_\alpha/\rho_t$ we get
\begin{multline}
g^{\alpha\beta}\phi_\alpha\phi_\beta\, \rho_t -2\phi_t
g^{\alpha\beta}\phi_\alpha \rho_\beta=\phi_t^2 \rho_t\big(
g^{\alpha\beta}\widehat{\phi}_\alpha\widehat{\phi}_\beta\, -2
g^{\alpha\beta}\widehat{\phi}_\alpha \widehat{\,\rho}_\beta\big)
=\phi_t^2 \rho_t\big(
g^{\alpha\beta}(\widehat{\phi}_\alpha-\widehat{\,\rho}_\alpha)
(\widehat{\phi}_\beta-\widehat{\,\rho}_\beta)\, -
g^{\alpha\beta}\widehat{\rho}_\alpha \widehat{\,\rho}_\beta\big)\\
=g^{ij}\big({\phi}_i-\widehat{\,\rho}_i \phi_t\big)
\big({\phi}_j-\widehat{\,\rho}_j\,\phi_t\big)\rho_t\, -
g^{\alpha\beta}{\rho}_\alpha{\,\rho}_\beta \phi_t^2/\rho_t
\end{multline}
Moreover
\beq
g^{\alpha\beta}\phi_\alpha\phi_\beta -2\phi_t
g^{\alpha 0}\phi_\alpha=-g^{00} \phi_t^2+g^{ij} \phi_i\phi_j
\eq
Hence
\begin{multline}
g^{\alpha\beta}\phi_\alpha\phi_\beta\, w_t -2\phi_t
g^{\alpha\beta}\phi_\alpha w_\beta
=\frac{\kappa\nup\varepsilon\ln{|1+t|}}{|\rho-2|^{1+\nup}}
w\Big(g^{ij}\big({\phi}_i-\widehat{\,\rho}_i \phi_t\big)
\big({\phi}_j-\widehat{\,\rho}_j\,\phi_t\big)\rho_t\, -
g^{\alpha\beta}{\rho}_\alpha{\,\rho}_\beta \phi_t^2/\rho_t\Big)\\
+\frac{\kappa\varepsilon}{(1+t)|\rho-2|^\nup} w
\Big(-g^{00} \phi_t^2+g^{ij} \phi_i\phi_j\Big)
\end{multline}
Since $|H|<1/2$ it also follows that
$$
\frac{1}{2} (\phi_t^2+\delta^{ij}\phi_i\phi_j) \leq
-g^{00}\phi_t^2+g^{ij}\phi_i\phi_j\leq
2(\phi_t^2+\delta^{ij}\phi_i\phi_j)
$$
Moreover; \beq \int_0^t \int \phi_t \Boxr_g \phi \, w\, dx d\tau
\leq \int_0^t\!\!\frac{\,\,c \varepsilon\!}{1\!+\tau}
 \int_{\Sigma_\tau} \!\!\! |\pa\phi|^2\,w\,  dx \, d\tau
 +\frac{1\!}{c \varepsilon\,}\!\int_0^t\!\!
 \int_{\Sigma_\tau}\!\! (1\!+\tau)\,
 |\Boxr_g \phi|^2\,w\,  dx \, d\tau.
\eq
\end{proof}

\section{Poincar\'{e} lemmas with weights}
We note that $\pa_r \rho=\pa_p \rho+\pa_q \rho=(1-H_{LL}/4)\pa_q \rho$,
since $L_2^\alpha\pa_\alpha \rho=0$, so
the estimate \eqref{eq:eikonal6} for $\pa_q \rho$ also hold for $\pa_r \rho$
with $c_2$ replaced by $2c_2$.
The following lemma was essentially proven in
\cite{A2}:
\begin{lemma}\label{poincare}
Suppose that $w$ is as in \eqref{eq:energyweight} with $\kappa>2c_2/\nup$, and that
with $\nup>0$ as in \eqref{eq:energyweight}
\beq
 |\pa_\rho \pa_r\rho|\leq \frac{2c_2 \varepsilon\ln{|1+t|}}{
 (1+|\rho|)^{1+\nup}}\, \pa_r \rho,\qquad 0<\pa_r \rho<\infty
 \eq
 Then for functions supported in $r\leq t+1$ we have
 \beq \int \Big(\frac{|\phi|}{1+|\rho|} \frac{\pa \rho}{\pa r}\Big)^2
\,w\, dx +\int \Big(\frac{|\phi|}{1+|r-t|}\Big)^2 \,w\, dx\leq
 32\int |\pa \phi|^2\, w\, dx
 \eq
\end{lemma}
\begin{proof} It suffices to prove the estimate for the first integral since the
second estimate is a special case of the first with $\rho=r-t$.
If we introduce polar coordinates $\rho=\rho(r,t,\omega)$ and change variables $r=r(\rho)$ for fixed $(t,\omega)$ we get
 \begin{multline}
\int_0^\infty \Big(\frac{|\phi|}{|\rho-2|} \frac{\pa \rho}{\pa r}\Big)^2 \,w\,
 r^2 d r =\int_{-\infty}^1 \Big(\frac{|\phi|}{|\rho-2|}\Big)^2 \frac{\pa \rho}{\pa r} \,w\,
 r^2 d \rho
 =\int_{-\infty}^1 |\phi|^2 \frac{\pa \rho}{\pa r} \,w\,
 r^2 \big(\frac{\pa}{\pa\rho} \frac{1}{|\rho-2|}\big) \, d \rho\\
 =-2\int_{-\infty}^1 \frac{\phi}{|\rho-2|}\frac{\pa\phi}{\pa\rho}\frac{\pa \rho}{\pa r} \,w\,
 r^2 d \rho-\int_{-\infty}^1 \Big(\frac{\phi}{|\rho-2|}\Big)^2|\rho-2| \frac{\pa}{\pa\rho} \Big(
 \frac{\pa \rho}{\pa r} \,w\, r^2\Big) d \rho
 \end{multline}
 Because of the conditions above
 \begin{multline}
\frac{\pa}{\pa\rho} \Big(
 \frac{\pa \rho}{\pa r} \,w\, r^2\Big)
 =\Big(\frac{\pa}{\pa\rho} \frac{\pa \rho}{\pa r}\Big)\,w\, r^2
+\frac{\pa \rho}{\pa r} \big((\pa_\rho w) r^2+  2rw \pa_\rho r\big)\\
\geq -\frac{c_2\nup \varepsilon\ln{|1+t|}}{
 (1+|\rho|)^{1+\nup}}\, \pa_r \rho\,\, w \, r^2
+\frac{\kappa\nup\varepsilon\ln{|1+t|}}{|\rho-2|^{1+\nup}} \,\pa_r \rho\,\, w\, r^2
+2r w\geq 0.
 \end{multline}
 Therefore
 \beq
\int_0^\infty \Big(\frac{|\phi|}{|\rho-2|} \frac{\pa \rho}{\pa r}\Big)^2 \,w\,
 r^2 d r\leq 2 \Big(\int_0^\infty \Big(\frac{|\phi|}{|\rho-2|} \frac{\pa \rho}{\pa r}\Big)^2 \,w\,
 r^2 d r\Big)^{1/2}
\Big(\int_0^\infty \Big(\frac{\pa\phi}{\pa \rho} \frac{\pa \rho}{\pa r}\Big)^2 \,w\,
 r^2 d r\Big)^{1/2}
 \eq
and it follows that
\beq
\int_0^\infty \Big(\frac{|\phi|}{|\rho-2|} \frac{\pa \rho}{\pa r}\Big)^2 \,w\,
 r^2 d r\leq 4
\int_0^\infty \Big(\frac{\pa\phi}{\pa \rho} \frac{\pa \rho}{\pa r}\Big)^2 \,w\,
 r^2 d r
=4
\int_0^\infty \Big(\frac{\pa\phi}{\pa r}\Big)^2 \,w\,
 r^2 d r
 \eq
and the lemma follows from also integrating over the angular variables.
 \end{proof}

\section{Energy estimates for the nonlinear problem}
 We will now show energy bounds assuming the strong decay estimates. Let
\beq
 E_{k,i}(t)=\sum_{|\bold{k}|\leq k,\, |I|\leq i }
 \int |\pa{\,}\pa^{\bold{k}}Z^I \phi|^2\, w\, dx,
\eq
where $w$ is as in Proposition \ref{curvedenergyest} with $\kappa=2c_2/\nup$ so
the conditions in Proposition \ref{curvedenergyest} and
Lemma \ref{poincare} hold if $\varepsilon>0$ is sufficiently small.
\begin{prop}\label{energybound} Let $N\geq 14$ and set $N^\prime=[N/2]+2$.
Suppose that $\phi$ is a solution of $\Boxr_{g(\phi)}\phi=0$ for $0\leq t<T$
such that $\phi(t,x)=0$ when $|x|\geq t+1$. Suppose also that
\begin{align}
|\pa \phi|& \leq \frac{c_1\varepsilon}{1+t}\label{eq:firstderdecay},\\
|\pa^2 \phi|& \leq \frac{c_2\varepsilon}{1+t}\Big|\frac{\pa \rho}{\pa q}\Big|\frac{1}{(1+|\rho|)^{1+\nup}},\qquad \nup>0,
\label{eq:secondderdecay}\\
|\phi|+|Z\phi|&\leq \frac{c_2\varepsilon}{1+t}\big( (1+|q|)+(1+t)^{c_2\varepsilon}\big)\label{eq:firstvectorfielddecay},\\
|\pa Z^I\phi| +(1+|q|)^{-1}|Z^I\phi| &
\leq \frac{c_{N^\prime} \varepsilon}{1+t}(1+t)^{c_{N^\prime}\varepsilon},\qquad
\text{for } |I|\leq N^\prime.\label{eq:manyvectorfielddecay}
\end{align}
Then there are constants $C_{k,i}$, depending only on
the constant above, such that for $0\leq t<T$;
 \beq\label{eq:energybound}
E_{k,i}(t)\leq 8 \sum_{\ell=0}^i E_{k+\ell,i-\ell}\,(0)
 (1+t)^{C_{k,i\,}\varepsilon},\qquad k+i\leq N.
 \eq
\end{prop}
\eqref{eq:energybound} will follow
 from \eqref{eq:nonlinearenergyest} below using induction and a Gronwall
 type of argument that we postpone.

\begin{prop}\label{energyest} Suppose that the assumptions in Proposition \ref{energybound} hold.
Then for $k+i\leq N$, $k,i\geq 0$;
 \beq\label{eq:nonlinearenergyest}
 E_{k,i}(t) \leq
 4E_{k,i}(0)
 +\int_0^t\frac{c_2^\prime \varepsilon}{1\!+\tau} E_{k,i}(\tau)\, d\tau
 +4\!\int_0^t \frac{c_{N^\prime}^2 \varepsilon\, (1\!+\tau)^{\,c_{N^\prime}\varepsilon}\!\!}{1\!+\tau}
 \big(E_{k+1,i-1}(\tau)+E_{k-1,i}(\tau)\big) \, d\tau,
 \eq
 where $E_{-1,n}=0$, $E_{m,-1}=0$.
 \end{prop}
 By Proposition \ref{curvedenergyest} with $\kappa=2c_2/\nup$
\beq\label{eq:energyest}
 E_{k,i}(t) \leq
 4E_{k,i}(0)
 +\int_0^t\frac{8(c_1\!+c_2) \varepsilon}{1\!+\tau} E_{k,i}(\tau)\, d\tau
 +4 \!\!\!\!\!\!\sum_{|\bold{k}|\leq k,\, |I|\leq i}\int_0^t \frac{1\!+ \tau}{c_2\varepsilon}
\int |\Boxr_g \pa^{\bold{k}} Z^I\phi|^2 \,w\,dx \, d\tau.
\eq
\subsection{Proof of \eqref{eq:nonlinearenergyest} in case $i=0$.}
 If $|\bold{k}|=0$ then $\Boxr_g\phi=0$ and \eqref{eq:nonlinearenergyest}
follows directly from \eqref{eq:energyest}, so we may assume that $1\leq |\bold{k}|\leq N$. If we use \eqref{eq:manycommutemixed} and
\eqref{eq:Hder2}, which holds since we assumed that $|\pa^{\bold{m}}\phi|\leq 1$ for $|\bold{m}|\leq N^\prime$, we get
 \beq
 |\Boxr_g \pa^\bold{k}\phi|\leq C\,|\pa \phi|\, \sum_{|\bold{n}|=|\bold{k}|+1
 \!\!\!\!\!\!\!\!\!\!\!} |\pa^\bold{n}\phi|\,+\sum_{|\bold{m}|+|\bold{n}|\leq |\bold{k}|+2, \,\,  1\leq |\bold{m}|\leq |\bold{k}|, \,\,  1\leq |\bold{n}|\leq |\bold{k}|
\!\!\!\!\!\!\!\!\!\!\!\!\!\!\!\!\!\!\!\!\!\!
\!\!\!\!\!\!\!\!\!\!\!\!\!\!\!\!\!\!\!\!\!\!\!\!\!\!\!\!\!\!\!\!\!}
 C \, |\pa^\bold{m}\phi|\, |\pa^\bold{n} \phi|,
 \eq
and hence by \eqref{eq:firstderdecay} and \eqref{eq:manyvectorfielddecay}
\beq\label{eq:boxkder}
 |\Boxr_g \pa^\bold{k}\phi|\leq \frac{c_1^\prime\varepsilon}{1+t} \sum_{|\bold{n}|=|\bold{k}|+1
 \!\!\!\!\!\!\!\!\!\!\!} |\pa^\bold{n}\phi|\,
 +\frac{c_{N^\prime} \varepsilon}{(1+t)^{1-c_{N^\prime}\varepsilon}}
  \sum_{1\leq |\bold{n}|\leq |\bold{k}|} |\pa^\bold{n} \phi|,
 \eq
since either $|\bold{m}|\leq N^\prime$ or $|\bold{n}|\leq N^\prime$,
in the second sum  in \eqref{eq:boxkder}.
\eqref{eq:nonlinearenergyest} in case $i=0$ follows from
\eqref{eq:energyest} using \eqref{eq:boxkder}.

\subsection{Proof of \eqref{eq:nonlinearenergyest} in case $k=0$.}
 By \eqref{eq:manycommute} and \eqref{eq:Hder2} using that we assumed that
 $|Z^J\phi|\leq 1$ for $|J|\leq [N/2]+2$, we have 3 types of terms:
 \beq
 |\Boxr_g Z^I\phi|\les
\les |Z^I\phi|\, |\pa^2\phi|
+\sum_{|J|\leq 1,\,\, |K|=|I|-1\!\!\!\!\!\!\!\!\!\!\!\!\!\!
\!\!\!\!\!\!\!\!\!\!\!\!\!\!\!\!\!\!\!\!\!} |Z^J \phi|\, |\pa^2 Z^K\phi|
+\!\!\sum_{|J|+|K|\leq |I|,\,\, |J|<|I|,\,\, |K|<|I|-1\!\!\!\!\!\!\!\!\!\!\!\!\!\!\!\!\!\!\!\!\!\!\!\!\!\!\!\!\!\!\!\!\!\!\!\!
\!\!\!\!\!\!\!\!\!\!\!\!\!\!} |Z^J \phi|\, |\pa^2 Z^K\phi|.
 \eq
By \eqref{eq:secondderdecay} and the Poincare lemma, Lemma \ref{poincare}
\beq
\int\Big(|Z^I\phi|\, |\pa^2\phi|\Big)^2 w \, dx\leq \int\Big(\frac{c_2^\prime\varepsilon}{1+t}\,
\Big|\frac{\pa\rho}{\pa q}\Big| \, \frac{|Z^I\phi|}{1+|\rho|}\Big)^2 w\, dx
\leq C\Big(\frac{c_2^\prime\varepsilon}{1+t}\Big)^2 \int |\pa Z^I\phi|^2 \, w\, dx.
\eq
By \eqref{eq:firstvectorfielddecay}
and \eqref{eq:tanZ} we have
\begin{multline}
\int \!\Big(\sum_{|J|\leq 1\!\!\!\!}|Z^J\phi| \sum_{|K|\leq |I|-1
\!\!\!\!\!\!\!\!\!\!\!\!\!\!} |\pa^2 Z^K \phi|\Big)^2 \! w\, dx
\leq \int\!\Big(\Big(\frac{c_2\varepsilon |q|}{1\!+\!t}\Big)^2\!\!\!
\sum_{|K|\leq |I|-1
\!\!\!\!\!\!\!\!\!\!}\!|\pa^2 Z^K \phi|^2+
\Big(\frac{c_2\varepsilon}{(1\!+\!t)^{1-c_2\varepsilon}}\Big)^2\!\!\!
\sum_{|K|\leq |I|-1
\!\!\!\!\!\!\!\!\!\!}\!|\pa^2 Z^K \phi|^2\Big)w\, dx\\
\leq \Big(\frac{c_2\varepsilon}{1+t}\Big)^2\int\sum_{|K|\leq |I|
\!\!\!\!\!\!\!}|\pa Z^K \phi|^2\, w\, dx+
\Big(\frac{c_2\varepsilon}{(1+t)^{1-c_2\varepsilon}}\Big)^2
\int \sum_{|K|\leq |I|-1
\!\!\!\!\!\!\!\!\!\!}|\pa^2 Z^K \phi|^2\, w\, dx
\end{multline}
The remain terms are easier to handle. Again by \eqref{eq:tanZ},
\eqref{eq:manyvectorfielddecay} and Lemma \ref{poincare}
\begin{multline}
\int\!\!\!\sum_{\begin{aligned}{}_{|J|+|K|\leq |I|,\,\, |J|<|I|,\,\, |K|<|I|-1}\end{aligned}
\!\!\!\!\!\!\!\!\!\!\!\!\!\!\!\!\!\!\!\!\!\!\!\!\!\!\!\!\!\!\!\!\!\!\!\!
\!\!\!\!\!\!\!\!\!\!\!\!\!\!\!\!\!\!\!\!\!\!\!\!\!\!\!\!}
\big(|Z^J \phi|\, |\pa^2 Z^K\phi|\,\big)^2 w\, dx\leq
\int\!\Big(\!\!\sum_{{\begin{aligned} {}_{|J|\leq |I|/2,\,\, |K|<|I|}\end{aligned}}\!\!\!\!\!\!\!\!\!\!\!\!\!\!\!\!\!\!\!\!\!\!\!\!\!\!\!\!\!\!\!}
\Big( \frac{|Z^J \phi|}{1+|q|}\, |\pa Z^K\phi|\,\Big)^2
+\sum_{\begin{aligned} {}_{|J|<|I|,\,\, |K|\leq |I|/2+1}\end{aligned} \!\!\!\!\!\!\!\!\!\!\!\!\!\!\!\!\!\!\!\!\!\!\!\!\!\!\!\!\!\!\!\!\!\!\!\!\!\!\!\!}
\Big( \frac{|Z^J \phi|}{1+|q|}\, |\pa Z^K\phi|\,\Big)^2 \Big) w\, dx\\
\leq
\Big(\frac{c_{N^\prime}\varepsilon}{(1+t)^{1-c_{N^\prime}\varepsilon}}\Big)^2
\int\sum_{|K|<|I|\!\!\!\!\!\!\!\!}|\pa Z^K\phi|^2 \, w\, dx
+\Big(\frac{c_{N^\prime}\varepsilon}{(1+t)^{1-c_{N^\prime}\varepsilon}}\Big)^2
\int \sum_{|J|<|I| \!\!\!\!\!\!\!}|\pa Z^J \phi|^2 \, w\, dx .
\end{multline}
Summing up we get
\begin{multline}
\int \!|\Boxr_g Z^I\phi|^2 w dx\\
\leq
\Big(\frac{c_2\varepsilon}{1\!+\!t}\Big)^2\!\!\!\int\!\!\sum_{|K|\leq |I|
\!\!\!\!\!\!\!}|\pa Z^K \phi|^2\, w\, dx+
\Big(\frac{c_2\varepsilon}{(1\!+\!t)^{1-c_2\varepsilon}}\Big)^2
\!\Big(\!\int\!\!\! \sum_{|K|\leq |I|-1
\!\!\!\!\!\!\!\!\!\!}\!|\pa^2\! Z^K \phi|^2\, w\, dx
+\int\!\!\! \sum_{|K|\leq |I|-1
\!\!\!\!\!\!\!\!\!\!}\!|\pa Z^K \phi|^2\, w\, dx\Big)
\end{multline}
\eqref{eq:nonlinearenergyest} in case $k=0$ follows from this using
\eqref{eq:energyest}.

\subsection{Proof of \eqref{eq:nonlinearenergyest} in case $k\geq 1$ and $i\geq 1$.}
Since $|\pa^\bold{m} Z^J\phi|\leq 1$ for $|\bold{m}|+|J|\leq N-5$, it follows from \eqref{eq:bigcommuteest} that
\beq
 |\Boxr_g \pa^{\bold{k}} Z^I\phi|\les
 \sum_{|\bold{n}|\leq |\bold{k}|,\, \, |J|+|K|\leq |I|,\,\,
|K|<|I|
 \!\!\!\!\!\!\!\!\!\!\!\!\!\!\!\! \!\!\!\!\!\!\!\!\!\!\!\!\!\!\!\!\!
 \!\!\!\!\!\!\!\!\!\!\!\!\!\!\!\!\!\!\!\!\!\!\!\!\!\!}   |Z^J \phi|\,\, |\pa^2 \pa^{\bold{n}} Z^K\phi|\,\,
+\sum_{|{\bold{m}}|+|\bold{n}|\leq |\bold{k}|,\,\,
 |J|+|K|\leq |I|
\!\!\!\!\!\!\!\!\!\!\!\!\!\!\!\!\!\!\!\!\!\!\!\!\!\!\!\!\!\!\!\!\!\!\!\!\!\!\!\!
\!\!\!\!\!\!\!\!\!\!\!\!\!\!\!\!\!\!}  |\pa \pa^{\bold{m}} Z^{J}\phi|\, |\pa \pa^{\bold{n}} Z^{K}\phi|,
\eq
for $|\bold{k}|+|I|\leq N$.
The terms in the first sum can be dealt with as in the case $k=0$ and
the terms in the second sum can be dealt with as in the case $i=0$.

\subsection{Proof of \eqref{eq:energybound} in case $k=i=0$} If $k=i=0$ then by \eqref{eq:energyest}
\beq
 E_{0,0}(t) \leq
 4E_{0,0}(0)
 +\int_0^t\frac{c_2^\prime \varepsilon}{1\!+\tau} E_{0,0}(\tau)\, d\tau
 \eq
and \eqref{eq:energybound} in case $i=k=0$ follows from this using
Lemma \ref{gronwall} below.

\subsection{Proof of \eqref{eq:energybound} in case $i=0$ and $k\geq n\geq 1$
 assuming \eqref{eq:energybound} in case $i=0$ for $k\leq n-1$}
By \eqref{eq:nonlinearenergyest} using \eqref{eq:energybound} for
$E_{k-1,0}$ we have
\begin{multline}
 E_{k,0}(t) \leq
 4E_{k,0}(0)
 +\int_0^t\frac{c_2^\prime \varepsilon}{1\!+\tau} E_{k,0}(\tau)\, d\tau
 +4\!\int_0^t \frac{c_{N^\prime}^2 \varepsilon\, (1\!+\tau)^{\,c_{N^\prime}\varepsilon}\!\!}{1\!+\tau}
 E_{k-1,0}(\tau) \, d\tau\\
\leq
4E_{k,0}(0)
 +\int_0^t\frac{c_2^\prime \varepsilon}{1\!+\tau} E_{k,0}(\tau)\, d\tau
 +32\!\int_0^t \frac{c_{N^\prime}^2 \varepsilon\, (1\!+\tau)^{\,(c_{N^\prime}+C_{k-1,0})\varepsilon}\!\!}{1\!+\tau}
 E_{k-1,0}(0) \, d\tau
\end{multline}
and again the estimate \eqref{eq:energybound} follows from
Lemma \ref{gronwall} below with $A=E_{k-1,0}(0)\leq E_{k,0}(0)$
and $B=c_2^\prime+32 c_{N^\prime}^2+ C_{N^\prime}+C_{k-1,0}$.

\subsection{Proof of \eqref{eq:energybound} in case $i=m\geq 1$ and $k=n\geq 1$
 assuming \eqref{eq:energybound} if $i=m$,  for $k\leq n-1$ and if
$i=m-1$ for all $k$, such that $i+k\leq N$.}
We will prove \eqref{eq:energybound} by induction in $i$ and for fixed $i$
induction in $k$. Since we have proven \eqref{eq:energybound} for $i=0$ and $k=0$ it suffices to prove \eqref{eq:energybound} in case $i=m\geq 1$ and $k=n\geq 1$
 assuming \eqref{eq:energybound} if $i=m$,  for $k\leq n-1$ and if
$i=m-1$ for all $k$, such that $i+k\leq N$.
By \eqref{eq:nonlinearenergyest} using \eqref{eq:energybound} for
$E_{k-1,i}$ and for $E_{k+1,i-1}$, we have
\begin{multline}
 E_{k,i}(t) \leq
 4E_{k,i}(0)
 +\int_0^t\frac{c_2^\prime \varepsilon}{1\!+\tau} E_{k,i}(\tau)\, d\tau
 +4\!\int_0^t \frac{c_{N^\prime}^2 \varepsilon\, (1\!+\tau)^{\,c_{N^\prime}\varepsilon}\!\!}{1\!+\tau}
 \big(E_{k+1,i-1}(\tau)+E_{k-1,i}(\tau)\big) \, d\tau\\
\leq 4E_{k,i}(0)
 +\int_0^t\frac{c_2^\prime \varepsilon}{1\!+\tau} E_{k,i}(\tau)\, d\tau
 +32\!\int_0^t \frac{c_{N^\prime}^2 \varepsilon\, (1\!+\tau)^{\,(c_{N^\prime}+C_{k+1,i-1}+C_{k-1,i})\varepsilon}\!\!}{1\!+\tau}
 \big(\widetilde{E}_{k+1,i-1}(0)+\widetilde{E}_{k-1,i}(0)\big) \, d\tau,
 \end{multline}
where $\widetilde{E}_{k,i}=\sum_{\ell=0}^i E_{k+\ell,i-\ell}$.
Using Lemma \ref{gronwall} with $A=\widetilde{E}_{k-1,i}(0)+\widetilde{E}_{k+1,i-1}$
and $B=c_2^\prime+32 c_{N^\prime}^2+ C_{N^\prime}+C_{k-1,i}+C_{k+1,i-1}$
we get with $C_{k,i}=2B$;
\beq
E_{k,i}(t)\leq \big( 4E_{k,i}(0)+\widetilde{E}_{k-1,i}(0)+\widetilde{E}_{k+1,i-1}(0)\big)
(1+t)^{C_{k,i} \varepsilon}\leq 8 \widetilde{E}_{k,i}(0)
(1+t)^{C_{k,i} \varepsilon}
\eq

We conclude by giving the Gronwall type of lemma used above:
\begin{lemma}\label{gronwall} Suppose that for some constants $A,B\geq 0$
\beq\label{eq:gronwallenergy}
E(t)\leq 4E(0)+\int_0^t \frac{B\varepsilon}{1+\tau}\Big( E(\tau)
+ A\,(1+\tau)^{B\varepsilon}\Big)\, d\tau.
\eq
Then
\beq\label{eq:gronwallest}
E(t)\leq (4E(0)+A)(1+t)^{2B\varepsilon}.
\eq
\end{lemma}
\begin{proof}
 If $G(t)$ denotes the integral in the right of \eqref{eq:gronwallenergy}
  then we have
 $$
 G^{\,\prime}(t)\leq\frac{B\ve}{1+t}G(t)+\frac{B\ve }{(1+t)^{1-B\ve}}A.
 $$
 If we Multiply with the integrating factor
 $$
 \frac{d}{dt}\Big( G(t) (1+t)^{-B\ve}\Big)
 \le \frac{B\ve}{1+t} A,
 $$
 and integrate we get
 \beq
G(t)(1+t)^{-B\ve}\leq G(0)+ B\ve A\ln{|1+t|}\leq G(0)+A (1+t)^{B\ve},
 \eq
and hence
\beq
E(t)\leq G(t)\leq G(0)(1+t)^{B\ve}+A(1+t)^{2B\ve} \leq (4E(0)+A )(1+t)^{2B\ve}.
\eq
\end{proof}

\section{Klainerman-Sobolev inequalities and $L^1-L^\infty$ estimates}\label{klainermansobolev}
First we state the Klainerman-Sobolev inequality:
\begin{prop} \label{klainerman} We have
 \beq
(1+t+|x|)(1+||t|-|x||)^{1/2}|\phi(t,x)|\leq C\sum_{|I|\leq 2}
{||Z^I \phi (t,\cdot)||_{L^2}}.
 \eq
\end{prop}
Next we state an inequality due to H\"ormander:
\begin{prop}\label{hormander}
 Suppose that $w(0,x)=\pa_t w(0,x)=0$. Then
 \beq
|w(t,x)|(1+t+|x|) \leq C\sum_{|I|\leq
2}\int_0^t\int \frac{|(Z^I \Box w)(\tau,y)|}{1+\tau+|y|}\, dy \, d\tau.
 \eq
\end{prop}
\begin{cor} \label{hormandercor}
Suppose that $\phi(0,x)=\pa_t \phi(0,x)=0$, when $|x|\geq 1$.
Then
 \beq
|\phi (t,x)|(1+t+|x|) \leq C\sum_{|I|\leq
2}\int_0^t\int \frac{|(Z^I \Box \phi)(\tau,y)|}{1+\tau+|y|}\, dy \, d\tau
+C\sum_{|I|\leq 2} \|\pa Z^I \phi (0,\cdot)\|_{L^2} .
 \eq
\end{cor}
\begin{proof} The inequality follows from writing $\phi=v+w$,
where $\Box w=\Box \phi$, $w(0,x)=\pa_t w(0,x)=0$, and
$\Box v=0$, $v(0,x)=\phi(0,x)$, $\pa_t v(0,x)=\pa_t \phi(0,x)$.
The inequality for $w$ follows from Proposition \ref{hormander}
and we will argue that the inequality for $v$ also follows from
Proposition \ref{hormander}. The inequality for $0\leq t\leq 1$ follows from
the usual Sobolev's lemma so it remains to prove it for $t\geq 1$. Let $\chi(t)$ be a smooth cutoff function
so that $\chi(t)=0$ when $t\leq 0$ and $\chi(t)=1$ when $t\geq 1$.
Then $\Box (\chi v)= \chi^{\prime\prime} v+2\chi^\prime v_t$
is supported in the set where $0\leq t\leq 1$ and $|x|\leq 2$ and it has vanishing initial data. It therefore follows Proposition \ref{hormander}
applied to $\chi v$ that for $t\geq 1$; $|v(t,x)|(1+t+|x|)=|\chi v (t,x)|(1+t+|x|)\leq C\|\Box (\chi v)\|_{L^1}
\leq C\sup_{0\leq t\leq 1} \|\pa v(t,\cdot)\|_{L^2}= C\|\pa v(0,\cdot)\|_{L^2}$.
\end{proof}

\section{The continuity argument}\label{proof}

Let $N\geq 14$ and set
\beq
E_N(t)=\sum_{|I|\leq N}\int |\pa Z^I\phi(t,x)|^2\, dx .
\eq
In view of local existence results it suffices to give a bound for
$E_N(t)$. We assume that initial data are so small that
\beq\label{eq:energydata}
E_N(0)\leq \varepsilon^2.
\eq
Fix $0<\delta<1$. We will argue by continuity. We assume the bound
\beq\label{eq:energyassumption}
E_N(t)\leq 16 N \varepsilon^2 (1+t)^\delta,
\eq
for $0\leq t\leq T$, which holds for $T=0$, and we will show that this bound
implies the same bound with $16$ replaced by $8$ if $\varepsilon$ is
sufficient small (independently of $T$).

Using Proposition \ref{klainerman} and \eqref{eq:energyassumption}
gives
\beq
|\pa Z^I\phi(t,x)|\leq \frac{C\varepsilon}{(1+t)^{1-\delta/2} (1+|t-r|)^{1/2}},
\qquad |I|\leq N .
\eq
Integrating this in the $t-r$ direction from $r-t=1$, where $\phi=0$ gives
\beq\label{eq:weakdecayenergy}
|Z^I\phi(t,x)|\leq \frac{C\varepsilon(1+|t-r|)^{1/2}}{(1+t)^{1-\delta/2}}
\leq \frac{C\varepsilon}{(1+t)^{(1-\delta)/2}}, \qquad
|I|\leq N-2 .
\eq
Since $\delta<1$ the weak decay estimate \eqref{eq:weakdecay} hold and hence
the decay estimates in Proposition \ref{strongdecay} as well
as the estimates for the solution of the approximate eikonal equation
in Proposition \ref{eikonalone} and the lemmas in the same section.
It therefore follows that the assumptions of Proposition \ref{energybound} hold
and that we therefore have
\beq\label{eq:finalenergyest}
E_N(t)=E_{0,N}(t)\leq 8\sum_{\ell=0}^N E_{\ell,i-\ell}(0)(1+t)^{C_{0,N}\varepsilon}\leq 8N E_N(0)(1+t)^{C_{0,N}\varepsilon}=8N\varepsilon^2(1+t)^{C_{0,N}\varepsilon},
\eq
since the family of vector fields $Z$ also contain the usual derivatives.
If $\varepsilon$ is so small that $C_{0,N}\varepsilon\leq \delta$,
then we get back the estimate \eqref{eq:energyassumption} with $16$ replaced
by $8$. This concludes the proof of  \eqref{eq:energyassumption}
and hence of Theorem \ref{mytheorem}. However, the proof above at most gives
the weak decay estimate \eqref{eq:weakdecay},
and hence the strong decay estimates in Proposition \ref{strongdecay}, with $\nu=1/2-c\varepsilon$, $c>0$. An additional argument using Corollary \ref{hormandercor} easily gives the weak decay estimate \eqref{eq:weakdecay} with $\mu=1-c\varepsilon$. In fact, since $\Boxr_g \phi=0$ we have using
\eqref{eq:weakdecayenergy};
\beq
|\Box Z^I\phi|\leq C\sum_{|J|+|K|\leq |I|} |Z^J\phi| \, |\pa^2 Z^K\phi |
\leq C\sum_{|J|+|K|\leq |I|+1, |J|\leq |I|} \frac{|Z^J\phi|}{1+|q|}
|\pa Z^K\phi|,
\eq
and hence using Lemma \ref{poincare} and H\"older's inequality
\beq
\sum_{|L|\leq 2} \int |Z^L \Box Z^I \phi(t,x)|\, dx
\leq C E_N(t) , \qquad\text{if}\quad |I|\leq N-3.
\eq
Therefore by Corollary \ref{hormandercor} we have for $|I|\leq N-3$
\beq
|Z^I\phi(t,x)|(1+t+|x|)\leq C\int_0^t \frac{E_N(\tau)}{1+\tau}\, d\tau+
CE_N(0)\leq C\varepsilon(1+t)^{C_{0,N}\varepsilon}.
\eq

\end{document}